\documentclass[14pt,a4paper]{article}
\usepackage[utf8]{inputenc}
\usepackage[OT1]{fontenc}

\usepackage{amsmath, amssymb, dsfont,amsthm}
\usepackage{enumitem}
\newcommand{\lcl}{cell-based {\m population control}}
\newcommand{\indi}[1]{\mathds{1}_{#1}}

\newcommand{\N}{\mathbb{N}}
\newcommand{\R}{\mathbb{R}}
\newcommand{\proba}{\mathbb{P}}
\newcommand{\F}{\mathcal{F}}

\newcommand{\nobj}{N_{obj}}
\newcommand{\wobj}{w_{obj}}
\newcommand{\wobjm}{w_{obj}^m}
\newcommand{\nmax}{N_{max}}
\newcommand{\target}{X_\mathbb{T}}
\newcommand{\targete}{{\mathcal T}_\epsilon}

\usepackage[linesnumbered]{algorithm2e}
\usepackage{algorithmic}
\makeatletter
\renewcommand{\@algocf@capt@plain}{above}
\makeatother

\usepackage{geometry}
\usepackage{graphicx}

\usepackage{float}

\usepackage{authblk}

\usepackage{xcolor}
\definecolor{darkgreen}{rgb}{0.1,.8,0.1}
\definecolor{darkred}{rgb}{0.8,.1,0.1}
\definecolor{darkblue}{rgb}{0.1,.1,0.8}
\definecolor{darkbluegreen}{rgb}{0.1,.6,0.8}

\newcommand{\m}{}
\newcommand{\md}{}
\newcommand{\modif}{}
\newcommand{\modiftwo}{}

\graphicspath{{Graphiques/Sur_vecteurs/},{Graphiques/Graphiques_Marshak/}}

\newtheorem{theorem}{Theorem}[section]

\newtheorem{lemma}[theorem]{Lemma}
\newtheorem{proposition}[theorem]{Proposition}
\newtheorem{remark}[theorem]{Remark}

\begin{document}
\date{May 12th, 2022}
\title{A cell-based {\m population control} of Monte Carlo particles
	{\m for the} 
	global variance reduction  for transport equations}
\author[,1]{Laetitia Laguzet\thanks{Corresponding author: \texttt{laetitia.laguzet@cea.fr}}}
\author[2]{Gabriel Turinici}

\affil[1]{CEA-DAM-DIF,F-91297 Arpajon, France}
\affil[2]{CEREMADE – UMR 7534, Université Paris Dauphine, PSL Research University, France}

\maketitle

\begin{abstract}
We present a population control method with sampling and regulation steps for Monte Carlo particles involved in the numerical simulation of a transport equation.
We recall in the first section the difficulties related to the variance reduction methods in the general framework of transport equations; we continue 
with a brief presentation of the mathematical tools invoked when  solving the radiative transport  equations and we focus on the 
importance of the emission and {\m control of existing} Monte Carlo particles.

The next part discusses several novel methods based on the \lcl \ method proposed in \cite{laguzet2020}. 
To this end, we analyze theoretically two types of splitting: one is conservative in energy (at the particle level) and the other is not.
Thanks to these results, a new algorithm is introduced that uses \lcl \  and a spatial distribution.
A numerical comparison of the different types of splitting is proposed in a simplified framework, then the various algorithms 
presented are compared against two benchmarks~:
the propagation of a Marshak wave
and
the propagation of two waves having different intensity and speed scales.
To carry out these last tests, we use the multi-physics code FCI2 \cite {dattolo01}.
\end{abstract} 

\noindent{\bf Keywords:} variance reduction ; Monte Carlo method ; sampling method ; photon transport ; radiative transfer.

\section{Introduction \label{sec:intro}}

The goal of this work is to propose a Monte Carlo sampling and {\m population control} method for the simulation of transport equations. To better situate the contribution,
we review the existing literature on the {\m population control and other} variance reduction techniques used for the transport equations together with the stochastic approaches 
employed for the resolution of radiative transfer equations.
Then the section \ref{sec:math_frame} gives a brief introduction to the Implicit Monte Carlo approach \cite{FLECK1971313} (used later for the numerical tests) while the section  \ref{sec:plan} contains the detailed overview of this paper.

\subsection{State of the art and motivation \label{sec:art}}

The radiative transport equations lead to the resolution of a transport equation coupled to an equation describing the evolution of the internal matter density. 
To solve such a non-linear transport problem two types of algorithms have been proposed: deterministic \cite{bailey2008piecewise,MOREL1996445} and stochastic.
 These methods can be further classified {\md as belonging to several approaches. In the first approach the weight of an emitted particle is treated as an energy} (e.g. the Fleck linearization \cite{FLECK1971313} {\md using some implicit treatment of the emitted term to linearize the transport equation}); in another approach (e.g., \cite{Nkaoua93}) {\md the inversion of a matrix is necessary to reconstruct the emitted term; other methods use} the diffusion limit \cite{Kaoua91, CLOUET2003139}.
 More generally Lapeyre \& al. \cite{lapeyre1998methodes} describe two probabilistic approaches 
to solve a transport equation (the section \ref{sec:math_frame} briefly describes one such method) in a non-stationary case.

The stochastic Monte Carlo methods have the drawback to 
display a slow convergence and many works focus on improving the error that  scales as $C/\sqrt{N}$; here $C$ 
represents the error in the Monte Carlo population and $N$ is the number of particles. The improvement takes the form of variance reduction methods, of which one can distinguish three main variants, see \cite{mcgrath1973techniques}~: the modified sampling method (with the population control), the use of analytical equivalence or other special techniques {\m(as \textit{Sequential Sampling} ; \textit{Orthonormal Functions} ; \textit{Adjoint Method} ; \textit{Transformations} or \textit{Conditional Monte Carlo})}.

The choice of a method is often case dependent. The modification of the weight of a particle (see section \ref{sec:math_frame} for details) 
by an \textit{well chosen} importance function is a crucial part of the variance reduction techniques (see \cite{kalos2008monte, lux1991, kalos1986monte}). On the other hand the adjoint operator has been used for a long time (see \cite{glasstone1970}) and motivated the development of several approaches to select a suitable  importance function depending on the situation
\cite{avatar, HAGHIGHAT200325, wagner1998, smith2005, mosher2009, wagner2009}.
A Russian-Roulette and Splitting method  (see \ref{sec:math_frame} for the definition), coupled to an importance function allows to conserve the weight 
in a window (depending on the geometry) and improves the computation
\cite{lux1991, booth1985, peplow2007,wagner2007,wagner2009,wagner2009bis, wagner2011,wagner2014,
	bal2018,turner1,turner2,LEGRADY2020107752} {\m and  \cite{doi:10.13182/NSE00-34}}.
For a more detailed description of the existing methods in the stationary case see for instance the first part of \cite{kulesza2018}.

The cited algorithms treat a stationary problem, using Russian Roulette (henceforth abbreviated RR) and Splitting (abbreviation S)
during the transition from a region of the domain to another (see \cite{juzaitis80, hoogenboom2004} for an analysis of these methods in this framework). On the contrary, very few works treat the {\m non-stationary} case and the variance reduction techniques have to be adapted to the situation under consideration (see \cite{Aa85, blanc2018,zimmerman1997monte}), or modified IMC approach with deterministic results \cite{urbatsch2012reduced,wollaber2009hybrid} {\md and  \cite{MBA, HOLO}}.
The presence of a source term implies the emission of particles; the computation cost has to be controlled, which is a concern not present in the stationary situation for the radiative transport problem (that is not the case for neutronic simulations).

{\modiftwo 
In the non stationary case such as the radiative transport problem, a parallel \textit{comb method} was proposed in \cite{osti_883456} (see \cite{kuilman2012} for an application of this method to a neutronic problem and \cite{wollaber2016,booth1996, gray1998} for additional details on previous non parallel proposals).}
Variance reduction being favorably impacted by uniform weight distribution,
the comb method tends thus to kill or split particles in order to obtain a fixed (target) number on a given area and fix the weights of the obtained particle population. 
{\md The main difference between our method and the comb method is that} the comb method allows to obtain exactly the target number of particle {\md by using a unique random number per cell} while our method allows to obtain this number only in average {\md because it uses a random number by particle} (see part \ref{sec:met_lcl}); that implies a different final weight {\md which is adjusted during the renormalization step. This approach removes the dependence on the distribution of the particles and allows to perform the algorithm independently of each particle (and is therefore compatible with parallel processing)}.
{\m Our algorithm \ref{algo:wobj_homo} is similar in spirit to the method in \cite{Urbatsh2008}: the distribution of the particles follows the energy in the cell. 
However,
the algorithm that we propose is not iterative and does not look for an exact number of particles to distribute; the convergence is obtained through stochastic limits and 
not during the research of a fixed point. 
Our approach is probabilistic, which allows, on average, to guarantee the number of particles distributed on the domain. Nevertheless, to guarantee the robustness of the algorithm, a security is imposed so that, whatever the simulation, the number of particles distributed on the domain is not higher than the number requested by the user: the sum of the objective numbers per cell, on the whole domain, does not exceed the number requested. 
Lastly, we propose a method to reach the target number of particles per cell (through the algorithm \ref{algo:meth_lcl}).
Nevertheless, the number of particles per cell after the regulation phase is always stochastic so that the total number finally obtained can exceed the desired number fixed by the algorithm for the same reasons as in the case of a constant objective number by cell (cf. algorithm \ref{algo:meth_lcl}). 
}

The previous studies show that there are some difficulties when parametrizing the RR and S methods. The \lcl \ method proposed here
and described in section \ref{sec:met_lcl} requires a single parameter (that is a proxy for the computation cost) and can be adapted to all types of unsteady computations 
(2D or 3D geometry etc...) that needs a control of the Monte Carlo population. This method, see \cite{laguzet2020}, has been developed 
in the framework of a multi-physics computation code, requiring a correct statistics on the whole domain: the \lcl \ method does not aim to 
improve the result in a direction or area but the whole statistic. This method is easy to implement, does not require any specific 
``know-how'', does not compute any importance function, no adjoint or any other, potentially costly, resolution. In addition it is independent of the 
geometry (be it 2D, 3D etc) and since it is also independent of the Monte Carlo estimate (such as indicated value, passage time, form functions, see \cite[page 72]{lapeyre1998methodes} for details), it does not augment the computation time for the particles involved. Another advantage of the method is that it is not
related to the transition kernel (in the unsteady state coupled to a different physics the computation of the solution between  $t$ and  $t+\Delta t$ 
depends on the transition kernel before $t$ only through the solution obtained); the method can also be used complementary to source bias methods (ideally that do not bias using the weights but using the position or direction). It is possible to relate this method with the adaptive weight windows of the stationary case: at the beginning of the time step, one limits the variance of the weights (which is transcribed by imposing a maximal value to {\m the objective weight (noted $\wobj$}), see part \ref{sec:met_lcl}).

The method proposed here is novel by the fact that it applies to the whole Monte Carlo particle population: the source emission is adapted to the present population which is itself regulated (the points 1 and 2 of the algorithm presented in part \ref{sec:math_frame}). Moreover, this approach does not affect the tracking phase and does not require any importance region. Of course, it can also be used when there is no source term.

\subsection{Mathematical framework \label{sec:math_frame}}

We illustrate the utilization of the \lcl \ method
in the framework of radiative transfer using the Implicit Monte Carlo approach, see \cite{FLECK1971313,gentile2009implicit}. To this end, we briefly describe the method and the resolution algorithm that has been used and refer the reader to \cite{FLECK1971313,wollaber2016} for additional details.

To simplify the presentation, the radiative transfer model is used in the gray case that is, when the radiative intensity is averaged in frequency. We consider a general radiative transfer model where $I$, the radiative intensity, depends on the time $t \in \mathbb{R}^+$, the position $x \in \mathbb{R}^d$ and the direction of propagation $\Omega \in \mathbb{R}^d$ (with $||\Omega||_d = 1$) 
through the following 
equation (see \cite{duderstadt1979transport} for a general presentation of transport equation and \cite{mihalas1999,pomraning2005equations} for a global description of the radiation hydrodynamics equations):
\begin{equation}
\frac{1}{c} \partial_t I(t,x,\Omega) + \Omega \cdot \nabla I(t,x,\Omega) = j(t,x,\Omega) - k(t,x,\Omega) I(t,x,\Omega),
\label{eq:transport_grise}
\end{equation}
where $c$ is the speed of light, $j$ the emission term of the matter and $k$ the absorption coefficient. 
This last equation is coupled with the equation describing the internal {\md energy} density of matter which can be written, when the matter (of temperature $T$) is supposed {\m stationary} 
 (as its volumic mass is independent of time):
\begin{equation}
\partial_t(\rho \epsilon) = \rho C_V \partial_t T = \int_{\mathbb{S}^2} kId\Omega - \int_{\mathbb{S}^2} j d\Omega, 
\end{equation}
with $\left. C_V = \frac{\partial \epsilon}{\partial T}\right|_{\rho \ cst} > 0$ being the caloric capacity at constant volume and $\mathbb{S}^2$ the unit sphere of $\mathbb{R}^3$ .

Neglecting the hydrodynamic motion of the matter
we have $j(t,x,\Omega) = k(t,x) \mathcal{B}(T)$ with 
$\mathcal{B}(T) = \frac{acT^4}{4 \pi}$ the gray Plank constant; the equation \eqref{eq:transport_grise} reads now:
\begin{equation}
\frac{1}{c} \partial_t I + \Omega \cdot \nabla I = k\left( \frac{acT^4}{4 \pi} - I \right).
\end{equation}
 
We apply the Implicit Monte Carlo (IMC) methodology introduced in \cite{FLECK1971313} that is, we solve with a Monte Carlo method the linear transport equation (see \cite{duderstadt1979transport} for a general presentation of the transport equations):
\begin{equation}
\frac{1}{c} \partial_t I + \Omega \cdot \nabla I + k^n I = f^n k^n B^n + (1-f^n)k^n \int_{\mathbb{S}^2} I(t,x,\Omega') d\Omega',
\label{eq:transport_fleck}
\end{equation}
where $f^n = \frac{1}{1+\beta^nck^n\Delta t}$ is the Fleck coefficient with $\beta^n = \frac{4 a (T^n)^3}{\rho C_V^n}$.
The quantities $k^n$, $B^n$, $\beta^n$ depend on the space variable but are now constants in the time interval $[t^n, t^n+\Delta t[$.
 
Integrating equation \eqref{eq:transport_fleck} over all directions $\Omega$ and then on a cell $m$ of volume $V_m$ and on the time interval $[t^n, t^n+\Delta t[$, we obtain 
that the energy to be emitted on the cell $m$ during the time interval is : 
\begin{equation}
S^n_m = V_m f^n_mk^n_mac(T^n_m)^4 \Delta t.
\label{eq:emission_volumique}
\end{equation}
Moreover, denoting $\{w_{ini,p}^n\}_{ p \in \{1, \dots, N^{ini}\} }$ the set of particles present at time  $t^n$ (that is, resulting from the previous iteration) and $\{w_{vol,p}\}_{ p \in \{1, \dots, N^{vol}\} }$ the set of particles emitted during the interval $[t^n, t^n+\Delta t[$ we can express the density of radiative energy $E^n_{r,m}$ of the cell  $m$ at time $t^n$ and the emitted energy $S^n_m$ as:
 \begin{equation}
  E^n_{r,m} = \frac{1}{V_m} \sum_{i=1}^{N^{ini}} w_{ini,p}^n, \ \ 
  S^n_m = \sum_{i=1}^{N^{vol}} w_{vol,p}^n.
  \label{eq:Er_Svol}
 \end{equation}

The term $\int_{t^n}^{t^n + \Delta t} \int_m \int_{\mathbb{S}^2} \Omega \cdot \nabla I d \Omega dx dt$ is susceptible to use the boundary conditions of the equation \eqref{eq:transport_fleck} which are represented by the surface particles  
$\{w_{surf,p}\}_{ p \in \{1, \dots, N^{surf}\} }$.

To solve the linear equation \eqref{eq:transport_fleck}, one can use a direct Monte Carlo method \cite[section 3.2]{lapeyre1998methodes}, on which the \lcl \ method can be adapted, irrespective of the transition kernel (here an uniform random variable) or functions that are involved (as $k^n$, $f^n$ or $B^n$).

We briefly present an adaptation of the general direct method  with time and space discretization, presented in \cite[section 3.4]{lapeyre1998methodes}, adapting it to a multi-physics code and previous notation in algorithm \ref{algo:MC_directe}.

\begin{algorithm}
\caption{Summary of the direct Monte Carlo method with time and space discretization. \label{algo:MC_directe}}

\noindent \textbf{Initialization} : for each cell, sample $N_{ic}$ particles from the probability law $I(0, x, \omega) dx d \omega$ ;

\For{\bf each time interval $[t^n, t^n+\Delta t[$ :} { 
  $ \star$ (optional) {\m Control} of the present population $\{w_{ini,p}\}_{ p \in \{1, \dots, N^{ini}\} }$ with the contraint \eqref{eq:Er_Svol} \;  
 
  $ \star$ Emission of $N_{vol}$ (and eventualy $N_{surf}$) particles to take into account the source term (and eventualy the boundary conditions) during the time interval with the contraint \eqref{eq:Er_Svol}
  ; 
  
  $ \star$ Tracking of the Monte Carlo particles : each particle moves according to the equation \eqref{eq:transport_fleck} and the weight of each particle is updated with the general formula: 
  $w_i(t+\Delta t) = w_i(t)\exp\left(\int_t^{t+\Delta t} f^n k^n(X_i(s)) ds\right)$ where $X_i(s)$ is the position of the $i^{th}$ particle at time $s$ ;
  
  $\star$ computation of all quantities required for the other physics (as $E^{n+1}_{r,m}$) and update (due to other physics) of the equation 
  \eqref{eq:transport_fleck} (in particular the Fleck coefficient $f^{n+1}$ and the functions $k^{n+1}$ and $B^{n+1}$).
}
\end{algorithm}

A particle, once emitted, is followed  until it experiences one of the following events: it goes out of the domain, it is absorbed by the weight rule (alone) or is 
suppressed during the control phase.
In the corresponding Monte Carlo simulation, as presented here, the weight is an exponential function which decays to zero but would never disappear completely.

The law of large numbers implies convergence when the number of particles is {\m large}, however the computer capacity is limited. Therefore the question of population 
management arises: how many particles are to be emitted to represent well the source term ($S_m^n$) and regulate the existing population. 

We focus on the new ``\lcl \ method'' used in the FCI2 code (see \cite{laguzet2020} for a detailed presentation of the method) that combines an emission controlled by the value of the source term with a RR and S strategy occurring at the source term emission time. 

\subsection{Scope and structure of paper \label{sec:plan}} 
The section \ref{sec:intro} is dedicated to a general presentation of the framework of this contribution.
We start in section \ref{sec:art} with a presentation of the state of the art of variance reduction methods which are essentially adapted to the stationary transport equations. 
We only found very few references concerning the instationary case, which shows the difficulty of the parametrization of the RR and S steps, both from the point of view of variance reduction techniques and the point of view of the management of the particle population.

The part \ref{sec:math_frame} illustrates this difficulty in the framework of radiative transfer, by presenting the algorithm resulting from the utilization of the 
Implicit Monte Carlo method, followed by a direct Monte Carlo method. We present the different particle populations that are used in the new {\m cell-based  population control} method studied in the part \ref{sec:met_lcl}. The section \ref{sec:local_originale} describes a simplified version of the method based on 
 \cite{laguzet2020} and allows to present the advantages of this method.
 In order to understand the pertinence of such an approach, the section \ref{sec:algo_proposes} analyzes the convergence of the algorithm with different variants.
Then the section  \ref{sec:lcl_wobj} proposes a variant of the \lcl \ method by inverting the paradigm~: instead of fixing a target number of particles by cell, one chooses a target weight. There methods are compared in section  \ref{sec:res_num}. First, we illustrate the convergence speed of the two algorithms with or without conservative
splitting  (section \ref{sec:vect_norme}) then we propose an application of the three methods for two benchmark cases: the propagation of a Marshak wave (section \ref{sec:transfert_rad}) and the propagation of two waves.

\section{Several algorithms based on the \lcl \  method \label{sec:met_lcl}} 

In this part we present the original \lcl \ method  in the section \ref{sec:local_originale} 
and we study two variants.
The first one (sections \ref{sec:algo_proposes}  and \ref{sec:rq_terme_source}  with and without term source respectively) allows to converge to the same distribution as the \lcl \ method with improved speed. The second variant  (section \ref{sec:lcl_wobj}) proposes a different paradigm~: 
obtain a target weight homogeneous among the cells instead of a target number of particles.

\subsection{Description and general properties of the original \lcl \ algorithm \label{sec:local_originale}}

In this section we present the \lcl \ method following \cite{laguzet2020} that is used for a single cell~: indeed, in this version, the algorithm is local 
at the cell level and independent of the neighboring cells. We will prove that this method guarantees, in average, the presence of {\m a target (objective) number} 
 $\nobj$ {\m of} particles {\m with} similar weight by cell. We refer to \cite[section 1.2]{laguzet2020} for a detailed presentation of the method that we expose here in a simplified setting.
In addition, we do not consider $N^{surf}$ (and therefore the boundary condition) because its role is similar to  $N^{vol}$ and $S_m^n$.

\noindent The procedure described in the algorithm \ref{algo:meth_lcl} is enforced at the beginning of an iteration and replaces the steps  3 and 4 of the algorithm \ref{algo:MC_directe}. 
We recall the notations used:
\begin{itemize}
    \item $N_m^{ini}$ : the number of particles in the cell $m$ issued from the previous iteration;
    \item $\{w_p^{ini} \}_{p \in \{ 1, \dots, N_m^{ini} \}}$ : weight of the current particle population;
    \item $E_{r,m}^n = \sum_{p= 1}^{N_m^{ini}} w^{ini}_p$ : sum of weights of the particles in the cell $m$;
    \item $S_m^n$ : energy to be emitted during the iteration $n$;
    \item $N_{m,obj}^n$ : number of target particles in the cell (user parameter or parameter depending on the total computational cost). This integer can depend on time and space ({\m we just use $\nobj$ if there is no ambiguous for the time and cell)}.
\end{itemize}

\begin{remark}
The treatment of the term involving the initial conditions is similar to $S_m^n$ and will not be detailed further.
\end{remark}

\noindent One of the advantages of this approach is to consider simultaneously two populations of particles:
\begin{enumerate}
    \item that resulting from the previous time steps $\{w_p^{ini} \}_{p \in \{ 1, \dots, N_m^{ini} \}}$, that will undergo a RR or S phase;
    \item the particles that convey the source term $S^n_m$ of the current iteration: $\{w_p^{vol} \}_{p \in \{ 1, \dots, N_m^{vol} \}}$ .
\end{enumerate}

The only hypothesis required by the \lcl \ method is the positivity of the terms $w_p^{ini} \ \forall p \in \{1, \cdots N^{ini} \}$, (thus of $E_{r,n}^m$) and $S_m^n$.
For notation convenience, we omit to mention explicitly in the sequel the time dependence of the quantities specific to the \lcl \ method. {\m Moreover, standard notation $\indi{}$ for the indicator function is used.}

\begin{algorithm}
	\caption{Description of the \lcl \ method \label{algo:meth_lcl}}

   Computation of $\wobjm$ : target weight to be reached by any particle of the cell $m$ :
$\wobjm \leftarrow \frac{E_{r,m}^n + S_m^n}{N_{m,obj}^n}$. 

\If(computation of the number of particles $N_m^{vol}$ to represent $S_m^n$ :){$S_m^n > 0$ }
{
	 $N_m^{vol} \leftarrow \max \left( 1, \left\lfloor\dfrac{S^n_m}{\wobjm}\right\rfloor \right)$ \;
	 $ w_m^{vol} \leftarrow \frac{S_m^n}{N_m^{vol}}.$ \;
     emission during the time step of $N_m^{vol}$	 particles of weight $ w_m^{vol}$ \;
}

\If(Russian Roulette and Splitting for each particle depending on its weight){$E_{r,m}^n > 0$}{
   \For{$p \in \{1, \dots, N_m^{ini}\}$}{
$I_p \leftarrow \left\lfloor\dfrac{w_p^{ini}}{\wobjm}\right\rfloor $ \;
$R_p \leftarrow \dfrac{w^{ini}_p}{\wobjm} - I_p$ \;
$u^p_{01} \sim \mathcal{U}(0,1)$ (uniform)\;
\eIf(Russian Roulette){$I_p = 0$}{

\eIf(the particle is killed){$R_p < u^p_{01}$}{$w^{ini}_p \leftarrow 0$}(the particle survives){$w^{ini}_p \leftarrow \wobjm$}

}
(Splitting){ 
	 $N^{split}_p \leftarrow I_p + \indi{u_{01}^p < R_p}$ \;
	 Conservative splitting : 
	 $ w^{ini}_p \leftarrow \frac{w^{ini}_p}{N^{split}_p} $ \;
	 Non conservative splitting : 
	 $ w^{ini}_p \leftarrow \wobjm.$ \;
	 The particle $p$ is duplicated $N^{split}_p - 1 $ times. 
}
} 
}
$c_m \leftarrow \frac{E_{r,m}^n + S_m^n}{\sum_{p=1}^{N_m^{ini}} w_p^{ini} + N_m^{vol} \cdot w_m^{vol} } $ \;
\For(final renormalisation on the present particles){$p \in \{1, \dots, N_m^{ini} \}$}{
$w_p^{ini} \leftarrow c_m \times w_p^{ini} $}
\For(final renormalisation on the emitted particles){$p \in \{1, \dots, N_m^{vol} \}$}{
	$w_p^{vol} \leftarrow c_m \times w_p^{vol} $}
Non-void correction: when $S^n_m=0$ and if after the previous RR and S phases the number of resulting particles is zero, 
	then the number of particles is set to $1$ with weight $E_{r,m}^n$.
\end{algorithm}

\begin{remark}
Only the conservative splitting exactly conserves the energy of the particles. The Russian Roulette phase, required for the {\m existing population control}, does not strictly conserves the energy (sum of weights of particles) which motivate the steps 27 and 30 of the algorithm \ref{algo:meth_lcl}.
The renormalization step will correct any energy mismatch introduced by the non-conservative splitting.
\end{remark}

\begin{lemma} \label{lemma:conservation}
The Russian Roulette and Splitting (conservative or not) phases of the algorithm \ref{algo:meth_lcl} conserve, in average, the energy.
\end{lemma}

\noindent {\bf Proof of lemma \ref{lemma:conservation}:} After the Russian Roulette and Splitting steps and before the renormalization, the weight $\tilde{w}^{ini}_p$ of the particle can be written {\m (using notations of algorithm \ref{algo:meth_lcl})}:
\begin{itemize}
	\item if conservative splitting: \begin{equation}
	\tilde{w}_p^{ini} = \indi{\{w_p^{ini} < \wobjm \}} \indi{\{u^p_{01} < \frac{w_p^{ini}}{\wobjm} \}} \wobjm + \indi{\{w_p^{ini} \geq \wobjm \}} \frac{w_p^{ini}}{N^{split}_p} .
	\end{equation}
    \item if  non conservative splitting : \begin{equation}
    \tilde{w}_p^{ini} = \indi{\{w_p^{ini} < \wobjm \}} \indi{\{u^p_{01} < \frac{w_p^{ini}}{\wobjm} \}} \wobjm + \indi{\{w_p^{ini} \geq \wobjm \}} \wobjm .
\end{equation}
\end{itemize}
 
 Noting that  $\mathbb{E}[N^{split}_p] = \frac{w_p^{ini}}{\wobjm}$ we obtain then that the energy of a particle is conserved in average:
 \begin{equation}
 \mathbb{E}[\indi{\{w_p^{ini} < \wobjm \}} \tilde{w}^{ini}_p + N^{split}_p\indi{\{w_p^{ini} \geq \wobjm \}} \tilde{w}_p^{ini}  ] = w_p^{ini}.
 \end{equation}
 
The next section deals with the properties of this algorithm in its two variants: conservative or non-conservative splitting and, in order to compare them, analyzes their convergence
when it is repeatedly applied on the same population of particles. Then, based on the properties of the \lcl \ method a stronger variant is introduced.

\begin{remark}
	Note that if we use the non conservative splitting, the expression of $\tilde{w}^{ini}_p$ can simply be written by: $\tilde{w}^{ini}_p = \wobjm \indi{N^{split}_p \geq1}$ and the energy resulting from the phase of RR and S on the particle $p$ is then $N^{split}_p \times \wobjm$. 
\end{remark}

\subsection{Theoretical properties of the \lcl \ method without source term}\label{sec:algo_proposes}

We present in this section some theoretical guidance concerning the performance of the \lcl \ method. We  investigate what is the limit distribution 
when the  \lcl \ algorithm is repeatedly applied on the same population of particles. This situation can correspond to a simulation regime where the 
time step imposed by the physics is small enough such that the trajectory phase does not influence greatly the weights of the particles.
We give formal arguments (in a particular setting) to show that the limit distribution is uniform and the number of particles is $\nobj$.

Consider the process $(N^l_m)_{l \in \mathbb{N}}$ such that  $N^0_m = N^{ini}_m $ is relative to a set of weights $\{w^l_i\}_{i \in \{1, \dots, N^l_m\}}$ and  $\{w^0_i\}_{i \in \{1, \dots, N^0_m \}}$ = $\{w^{ini}_p\}_{p \in \{1, \dots, N^{ini}_m\}}$. 
The weights $\{w^{l+1}_i\}_{i \in \{1, \dots,N^{l+1}_m\}}$ are obtained by applying the algorithm \ref{algo:meth_lcl}  to the set of weights 
$\{w^l_i\}_{i \in N^l_m}$ which are not issued from the generation of the source term. 
The process  $(N^l_m)_{l \in \mathbb{N}}$ describes  the number of weights considered after $l$ iterations i.e., phases of RR and S
on the population of weights  $\{w^{ini}_p\}_{p \in \{1, \dots, N^{ini}_m\}}$ with the renormalization step (step at line 31 in algorithm \ref{algo:meth_lcl}).

\subsubsection{\m The convergence of the non-conservative splitting method}

We analyze the convergence of the algorithm \ref{algo:meth_lcl} presented in section \ref{sec:met_lcl} and more precisely the non-conservative splitting 
when $S^n_m = 0 $ (thus $N^{vol}_m=0$)
and, to simplify the presentation $E_{r,m}^n= 1$. Let us fix a cell. 

The algorithm
\ref{algo:meth_lcl} is local relative to each cell {\m and is applied at the beginning of a time interval}, thus for simplicity we omit the {\m indexes  $m$  and $n$} in this part.

\begin{proposition}  \label{prop:cv_nobj}
Denote by $N^l$ the number of the particles in the cell after the $l$-th iteration of the algorithm.  Then 
\begin{equation}
\lim\limits_{l \rightarrow +\infty}  \proba( N^l = \nobj) = 1.
\end{equation}
More precisely, there exists $\lambda > 0$ such that
\begin{equation}
\proba( N^l \neq \nobj) \le e^{- l \cdot \lambda },
\end{equation}
i.e., the convergence is exponential.
\end{proposition}

\begin{proof}
	 After $l > 1$ iterations of the algorithm, there are $N^l$ particles present of equal weight $E_{r,m}^n/N^l= 1/N^l$. Note that this common weight is not necessarily $\wobj = 1/\nobj$ due to the renormalization step. The number of particles that are
	 obtained at iteration $l+1$ is given by:
\begin{equation}
  N^{l+1} = \max(1, Y^l) \  \text{where }
 Y^l= N^{l} \times \left\lfloor\dfrac{\nobj}{N^{l}}\right\rfloor 
 + \mathcal{B}\left(N^{l}, \frac{\nobj}{N^{l}} - \left\lfloor\dfrac{\nobj}{N^{l}}\right\rfloor \right). 
\label{eq:defyl}
\end{equation}
\noindent 
We used the standard notations~:

- $\mathcal{B}(n,p)$ designates a binomial variable of parameters $n$ and $p$ independent of any previous other variable; 

- the floor function $\lfloor\cdot \rfloor$ designates the largest integer smaller than (or equal to) the argument.

The previous formula can also be written in the following form:
\begin{equation} 
N^{l+1} = \max\left(1, \nobj 
+ \sum_{y=1}^{\infty} \indi{N^{l} = y} \left[ \mathcal{B}\left( y,\epsilon_y \right) - y \epsilon_y\right] \right),
\label{eq:eqnl2}
\end{equation}
with 
\begin{equation}
\epsilon_y =  \frac{\nobj}{y} - \left\lfloor\dfrac{\nobj}{y}\right\rfloor. 
\end{equation}

We recall that the binomial variables appearing in \eqref{eq:eqnl2} are independent of the sigma algebra $\F_l = \sigma\{N^k, k\le l \}$ (the smaller sigma algebra making $N^0, ..., N^l$ measurable). 

Denote now $\max(N^k)$ to be the maximum value of the random variable $N^k$ (this value is finite because the random variable is discrete). 
We show next that $N^l$ are all bounded by $M_0 = \max(N^0, 2 \nobj)$. 

When 
$y \le \nobj$  all values taken by the random variable $\mathcal{B}\left( y,\epsilon_y \right)$ are smaller or equal to $y$, thus 
$\mathcal{B}\left( y,\epsilon_y \right) - y \epsilon_y \le y \le \nobj$. 

When $y \ge \nobj$  then $\mathcal{B}\left( y,\epsilon_y \right) - y \epsilon_y \le y - y\epsilon_y$. But in this case 
 $\epsilon_y =  \frac{\nobj}{y} $ thus   $y - y\epsilon_y  = y - \nobj$. 
Using \eqref{eq:eqnl2} we can write:
\begin{eqnarray} 
& \ & 
Y^{l} = 
\nobj + \sum_{y=1}^{\infty} \indi{N^{l} = y} \left[ \mathcal{B}\left( y,\epsilon_y \right) - y \epsilon_y\right] 
\\ \nonumber &  \ & 
\le  \nobj + 
\sum_{1 \le y \le \nobj} \indi{N^{l} = y} \cdot \nobj + \sum_{y > \nobj} \indi{N^{l} = y} \cdot (y- \nobj) 
\\ \nonumber &  \ & 
\le \indi{N^{l} \le \nobj} \cdot 2 \nobj + \sum_{y > \nobj} \indi{N^{l} = y} \cdot y
\\ \nonumber &  \ & 
\le  \indi{N^{l} \le \nobj} \cdot 2 \nobj +  \indi{N^{l}  > \nobj} \cdot \max(N^{l}) 
\le  \max(N^0, N^1, ..., N^l,2 \nobj).
\label{eq:eqnl3}
\end{eqnarray}
By recurrence we obtain $N^{l+1} \le M_0$.

For $a,b \le M_0$ we introduce the notation $p_{a\to b} := \proba(N^{k+1} = b | N^k=a)$ and remark that the quantity is independent of $k$. In fact $N^l$ 
are realizations of a Markov chain and $p_{a \to b}$ are the transition probabilities of this process. In particular note that $p_{a \to \nobj} > 0$ for any $a$ 
because the variable $Y^l$ takes all values in the interval $[\nobj -N^l, \nobj + N^l(1-\epsilon_{N^l})]$ which contains $\nobj$. Moreover 
$p_{\nobj \to \nobj} =1$. Denote also $p_- = \min \{ p_{a \to \nobj} ; a \neq \nobj \}$. Note that $p_- > 0$.
Thus

\begin{eqnarray} & \ & 
\proba(N^{l+1} = \nobj) =  \sum_{a \le M_0} \proba(N^{l+1} = \nobj | N^{l} = a )  \proba(N^{l} = a)
\\ \nonumber &  \ & 
= \proba(N^{l} = \nobj) + \sum_{a \neq \nobj}  p_{a \to \nobj}  \cdot \proba(N^{l} = a)
\\ \nonumber &  \ & 
\ge \proba(N^{l} = \nobj) + \sum_{a \neq \nobj}  p_-  \cdot \proba(N^{l} = a)
 \ge \proba(N^{l} = \nobj) +  p_- \cdot (1- \proba(N^{l} = \nobj)).
  \end{eqnarray}
We obtain thus 
\begin{equation}
\proba(N^{l+1} = \nobj)  \ge p_- + (1-p_-) \cdot \proba(N^{l} = \nobj),
\end{equation}
which can also be written
\begin{equation}
\proba(N^{l+1} \neq \nobj)  \le  (1-p_-) \cdot \proba(N^{l} \neq \nobj),
\end{equation}
and the conclusion follows with $\lambda = -\log(1-p_-)$.
\end{proof}

{\m

\subsubsection{ \m The convergence of the conservative method}

We continue with the analysis of the convergence of the conservative version of the  algorithm \ref{algo:meth_lcl} presented in section \ref{sec:met_lcl} when $S^n_m = 0 $ (thus $N^{vol}_m=0$)
and, to simplify the presentation, $E_{r,m}^n= 1$. Let us fix a cell. 
{\m For the same reasons as above, in} 
order to simplify the 
presentation we will omit the index $m$ of the cell{\m, the index $n$ of the time interval and} only keep the iteration count $l$ as index.
Since the cell is fixed and there are no source terms, the total mass, $E_{r,m}^n= 1$ is fixed~; $\nobj$ and $\wobj=E_{r,m}^n/\nobj = 1/\nobj$ are also fixed.

Denote from now on by $X_l$ the state of the algorithm after $l$ iteration of RR + S plus the renormalization steps.
Note that, in full generality, $X_l$ can be described as a set of (unknown number, noted  $N^l$, of) positive weights ({\m noted $X_l^1, \dots, X_l^{N^l}$}) that sum up to $E_{r,m}^n$, i.e. a member of $\ell^1$, the ensemble of (absolutely) summable sequences. 
Denote also by $\target$ the uniform distribution with $\nobj$ particles all of weight $\wobj$: 
$\target= (\wobj,...,\wobj)\in \R^{\nobj}$. With these preliminaries we can prove the following result.

\begin{proposition} \label{prop:cv_conservative}
	Suppose $N^0$ is finite. Then
	the sequence of random variables $(X_l)_{l\ge0}$ converges almost everywhere in $\ell^1$ to the target distribution $\target$ i.e.
	\begin{equation}
	\proba[ \lim_{l\to \infty} X_l = \target] = 1.
	\label{eq:convergence_conservative}
	\end{equation}
\end{proposition}

\begin{proof}

{\m The proof requires many technical arguments, we present here only a sketch and refer for the full proof to appendix \ref{sec:appendix_proof_nc}.

$\star$ First, we invoke lemma \ref{lemma:returntonobj} that ensures that, after a finite time $\tau$,  $X_l$ will have at most $\nobj$ non-null particles i.e. there exists $c_1>0$ such that:
\begin{equation}
\forall L \ge 0: \proba[ N_0 > \nobj, N_1 > \nobj,..., N_L > \nobj] \le e^{c_1 (L-1)}.
\end{equation}

$\star$ Then, we prove (lemma \ref{lemma:finitedimension}) that, for $l \geq \tau$, $X_l$ will have at most $\nmax= 6 \nobj$ non-null particles forever: for any $l\ge \tau$, $N^l \le 6 \nobj$ (with certainty).

$\star$ After that, we introduce the set 
(a specific neighborhood of the target distribution)~:
$$\targete = \left \{ (w_1,...,w_{\nobj}) \in \R^{\nobj} \left|
\sum_{k=1}^{\nobj} w_k= \nobj \wobj,
\sum_{k=1}^{\nobj} |w_k-\wobj| \le \epsilon \wobj \right.
\right\}$$
and denote $\tau_\epsilon$ the first iteration step $l$ when $X_l$ enters $\targete$.
In lemma \ref{lemma:convergencetargete} we show that there exists $\epsilon_0 >0$ such that 
for any $\epsilon \in[0, \epsilon_0]$ there exists $L_\epsilon \in \N$ and $c_\epsilon < 1 $ with:
\begin{equation}
\forall l \ge \tau: \ 
\proba[ 
X_{l+1} \notin \targete, ...,
X_{l+L_\epsilon} \notin \targete | X_{l} \notin \targete ] \le c_\epsilon.
\end{equation}
Thus $\proba[\tau_\epsilon \ge k]$ is exponentially decreasing when $k\to \infty$ and in particular the
stopping time $\tau_\epsilon$ is finite almost everywhere ($\proba[\tau_\epsilon < \infty]=1$) i.e., after a finite time $\tau_\epsilon$ the state $X_{\tau_\epsilon}$ will be in $\targete$.

$\star$ We conclude using lemma \ref{lemma:stability}:
there exists some constant $c_3 >0$ only depending on $\nobj$ such that
\begin{equation}
\proba[\lim_{l \to \infty} X_l = \target | \tau_\epsilon < \infty ] \ge 1-c_3 \epsilon.
\end{equation}
As this lemma is true for any $\epsilon$ (small enough) then we have convergence everywhere.

}

\end{proof}
\begin{remark}
	As mentioned before, in full rigor the convergence should be understood as the convergence of
	$\ell^1$ sequences, but in fact, as it will be seen later, it boils down to the usual convergence of
	$\nmax$-tuples in $\R^{\nmax}$. Also note that the result does not require any hypothesis on the number $N_0$ of particles present in the initial datum $X_0$ (in particular it can be as large as required, as long as it is finite; moreover, slight modifications of the proof allow to treat the situation of a infinitely countable set of non-null initial weights).
\end{remark}

}

\subsection{Theoretical properties of the \lcl \ method with source term}\label{sec:rq_terme_source}

In this part, we investigate the case when the source term cannot be neglected. In general the source term is time-dependent and therefore there is no
constant target distribution of the particles (in terms of number or weight). In order to give some theoretical guidance we will show that the 
RR+S (plus renormalization) part of the procedure is not affected by the choice of a objective weight $\wobjm$ depending on the source too (besides the previous energy);
under appropriate hypotheses we show that the RR+S part does {\m help converge the number of particles to} $\nobj$. 
We consider thus that the algorithm is repeated several times on the particles of the previous time step.
This is not the general case (in the general case one would apply the algorithm to all the particles) but supports the idea that by treating the emission simultaneously
with the {\m control of existing population}, 
the  \lcl \ method ensures a convergence towards $\nobj$ particles even with the addition of a source term. 

We reuse the notations from the previous part by distinguishing two populations, the initial one and the emitted one:
 $ N^l = N_m^{vol} + N_l^{ini} $ i.e., 
 {\m $N_{l+1}^{ini}$} 
 is the number of ``initial'' particles in the cell 
after {\m an} iteration of the RR + S part and renormalization of the 
$N_l^{ini}$ 
particles.

\begin{proposition} Suppose $\dfrac{S_m^n \nobj}{E_{r,m}^n+S_m^n} \in \N$ and  $N_m^{vol} >1$.
Consider the procedure that acts iteratively on two populations of particles as follows:
	
- 	on the $N_m^{vol}$ emitted particles defined at the line ``computation of the number of particles $N_m^{vol}$ to represent $S_m^n$'' of algorithm \ref{algo:meth_lcl} : each is assigned weight  $ w_m^{vol} = \frac{S_m^n}{N_m^{vol}}$ without further modifications;

-  on the non-volumic particles : one computes {\m $\wobj^m$} as detailed at line ``Computation of $\wobjm$'' of algorithm \ref{algo:meth_lcl}, runs the 
 the  RR+S steps on the $N_l^{ini}$ particles, and then renormalizes them to obtain total mass  $E_{r,m}^n+S_m^n$.
 
Denote $N^l$ the number of particles at the iteration $l$. Then:
	\begin{equation}
	\lim\limits_{l \rightarrow +\infty}  \proba[ N^l = \nobj ] = 1.
	\end{equation}
	More precisely, there exists $\lambda > 0$ such that
	\begin{equation}
	\proba( N^l \neq \nobj) \le e^{- l \cdot \lambda },
	\end{equation}
	i.e., the convergence is exponential.
\label{prop:cv_nobj2}
\end{proposition}

\begin{remark} The condition  $\dfrac{S_m^n \nobj}{E_{r,m}^n+S_m^n} \in \N$ is required if one hopes to find convergence to some uniform distribution while respecting the 
	conservation of mass. Indeed, if both the volumic and non-volumic particles have the same target weight  and there are $\nobj$ particles in all, then 
	$E_{r,m}^n$ and $S_m^n$ are both integer multiples of the target weight and 	relation 	is satisfied. On the other hand, the specificity of this proof, compared to 
	proposition \ref{prop:cv_nobj} is  that  ${\m \wobj^m}$ depends on the source term ${S_m^n}$ too.
\end{remark}	
	
\begin{proof} 
	After $l > 1$ iterations of the algorithm, there are $N_l^{ini}$ particles present of equal weight $\wobj/N^l$. Note that:
\begin{equation}
\begin{split}
N^{l+1}& = N^{vol} + N_{l+1}^{ini} , \text{where } N_{l+1}^{ini} = \max(1, Y^l) \text{ and }\\
& Y^l= N^{ini}_l \times \left\lfloor\dfrac{E_{r,m}^n}{E_{r,m}^n+S_m^n}\times \dfrac{\nobj}{N^{ini}_l}\right\rfloor 
+ \mathcal{B}\left(N^{ini}_l, \dfrac{E_{r,m}^n}{E_{r,m}^n+S_m^n}\times \dfrac{\nobj}{N^{ini}_l} - \left\lfloor\dfrac{E_{r,m}^n}{E_{r,m}^n+S_m^n}\times \dfrac{\nobj}{N^{ini}_l}\right\rfloor \right). 
\end{split}
\end{equation}

On the other hand, at the present iteration $l$  the number of emitted particles is $N^{vol}= \left\lfloor\dfrac{S_m^n \nobj}{E_{r,m}^n+S_m^n}\right\rfloor = \dfrac{S_m^n \nobj}{E_{r,m}^n+S_m^n}$, the last equality being true by hypothesis. 
Denote $\nobj^{ini} = \dfrac{E_{r,m}^n \nobj}{E_{r,m}^n+S_m^n} = \nobj - N^{vol}$ which is also an integer for the same reason. Then one can write:
\begin{equation}
N^{l+1} = N^{vol} + N_{l+1}^{ini} , \text{where } N_{l+1}^{ini} = \max(1, Y^l) \text{ and }
 Y^l= N^{ini}_l \times \left\lfloor \dfrac{\nobj^{ini} }{N^{ini}_l}\right\rfloor 
+ \mathcal{B}\left(N^{ini}_l, \dfrac{\nobj^{ini} }{N^{ini}_l} - \left\lfloor \dfrac{\nobj^{ini} }{N^{ini}_l}\right\rfloor \right). 
\end{equation}
But this formula is the same as \eqref{eq:defyl} with $\nobj$ replaced by $\nobj^{ini}$. Starting from this point the same proof {\m as in proposition \ref{prop:cv_nobj}} applies.
\end{proof}

\begin{remark} 
{\m Similarly, the proof for the conservative splitting is a combination of the above arguments and  proposition \ref{prop:cv_conservative}.}
\end{remark}

\begin{remark}
	When $\nobj \to \infty$ then $w_m^{vol} \to {\m \wobj^m}$. Indeed, $|w_m^{vol} - {\m \wobj^m}|= ({\m \wobj^m})^2/|(S_m^n-{\m \wobj^m})| \to 0 $ because ${\m \wobj^m} \to 0 $ 
	when $\nobj \to \infty$.    
\end{remark}

\subsection{Additional procedures based on the \lcl \ method \label{sec:lcl_wobj}}

The previous section allowed to show that the \lcl \ algorithm with conservative or non-conservative splitting will orient towards an uniform  distribution of 
the weights of the particles, locally in each cell. We already gave some details on the spatial and temporal dependence of $N_{m,obj}^n$ in section \ref{sec:met_lcl} 
but a too empirical 
parametrization of  $\nobj$ could impact the other variance reduction methods.

However, obtaining uniform weights during simulation could allow a global variance reduction; we propose thus to use the parametrization of $\nobj$ 
to obtain  objective weights by cell $\wobj$ as homogeneous as possible.

To this end, the total number of particles specified at the start of an iteration can be determined in several ways:
\begin{enumerate}
	\item given by a single user parameter, in order to best control the simulation cost;
	\item be derived from a user parameter as above. In this case, one has to specify the number of cells to be treated in Monte Carlo (according to the conditions dictated by the simulation) and multiply by $\nobj$.
\end{enumerate}

If the RR and S phases are global on the domain they can induce {\m a nonphysical phenomenon of energy transport} like a  teleportation error \cite{osti_1212834} {\m for the emitted term}. To limit it, 
 we keep a cell-based objective weight and only change the determination of the  $\nobj$ as described in the algorithm \ref{algo:wobj_homo}. 
Let $\mathcal{M}$ be the set of cells used in the simulation and $c_\mathcal{M}$ the number of cells.
We obtain the algorithm \ref{algo:wobj_homo} described below.
	 
\begin{algorithm}
	\caption{Modification of {\m the \lcl \ }method to obtain homogeneous weights {\m for all particles in mesh $m \in \mathcal{M}$}. \label{algo:wobj_homo}}
	
	{\noindent \textbf{Compute total particles number} :
	
		 case 1 : read the user parameter $N_{total}^n$ and set $N_{share}^n = N_{total}^n - 2 c_{\mathcal{M}}$ (see remark \ref{rq:Ntotal});
		 
		 case 2 : read the parameter $\nobj^{usr}$ to compute $N_{share}^n = (\nobj^{usr} -2) \times c_{\mathcal{M}} $ and $N_{total}^n = \nobj^{usr} \times c_{\mathcal{M}} $ (see remark \ref{rq:Ntotal})
}

	\noindent \textbf{Compute total energy} $E_{r,tot}^n = \sum_{m \in \mathcal{M}} E_{r,m}^n + S_m^n$ \;
	
	\For{$m \in \mathcal{M}$ }{
		\If{$E_{r,m}^n + S_m^n > 0$}{
		$\star$	$u^m_{01} \sim \mathcal{U}(0,1)$
			
		$\star$ $R_m = \frac{E_{r,m}^n - S_m^n}{E_{r,tot}^n} - \left\lfloor\dfrac{E_{r,m}^n + S_m^n}{E_{r,tot}^n}\right\rfloor$ 
		
		$\star$ $N_{m,obj}^n = \max \left[ \left\lfloor\dfrac{E_{r,m}^n + S_m^n}{E_{r,tot}^n}\right\rfloor N_{{\m share}}^n + \indi{u_{01}^m < R_m  } ; \indi{E_{r,m}^n >0} + \indi{S_m^n>0} \right]  $ ; \label{linealgo:linemax}
		
		$\star$ Apply algorithm \ref{algo:meth_lcl} with $\nobj = N_{obj,m}^n$.
	    }
	}
\end{algorithm}

\begin{remark} The goal of this method is to obtain homogeneous weights on the domain; the non conservative splitting is thus preferred.
\end{remark}

\begin{remark} \label{rq:Ntotal}
	Without using the ``max'' operation in line \ref{linealgo:linemax}, the algorithm \ref{algo:wobj_homo} can lead to a pathological situation as follows: 
	even if 
	$E_{r,m}^n >0$ or $S_m^n>0$, 	
	the number of particles obtained with the algorithm 
	can reach $N_{m,obj}^n = 0$ which may be unacceptable because the method must conserve the energy in each cell.
	That is why, the ``max'' operation is a fail-safe to make sure to not exceed $N^{total}$; it is set with the most pessimistic (and in fact impossible) alternative: 
	all cells are pathological. 
\end{remark}

\noindent With this method, previous properties  such as lemma \ref{lemma:conservation} are preserved. We can state:

\begin{lemma} \label{lemma:conservationNtotal}
The average total number of particles for iteration $n$ is less than $N_{total}^n$. Moreover, if $N_{m,obj}^n  = \left\lfloor\dfrac{E_{r,m}^n + S_m^n}{E_{r,tot}^n}\right\rfloor N_{{\m share}}^n + \indi{u_{01}^m < R_m  } > 0 \ \forall m \in \mathcal{M}$ then the average total number of particles for iteration $n$ is  $N_{share}^n$.
\end{lemma}

\begin{proof}
Let us denote by $\tilde{N}_m^n $ the random variable describing the number of particles in the cell $m$ with $N_{obj,m}^n \ge \indi{E_{r,m}^n >0} + \indi{S_m^n>0}$ after the action of the algorithm \ref{algo:meth_lcl}. Then: 
\begin{align}
\mathbb{E}[\tilde{N}_m^n|N_{obj,m}^n] &= N_m^{vol} + \mathbb{E}\left[ \sum_{p = 1}^{N_m^{ini}} \indi{\{w_p^{ini} < \wobjm \}} \indi{\{u^p_{01} < \frac{w_p^{ini}}{\wobjm} \}} + N^{split}_p\indi{\{w_p^{ini} \geq \wobjm \}}   \right] \\
&= \frac{S_m^n}{S_m^n+E_{r,m}^n} N_{obj,m}^n 
+ \sum_{p = 1}^{N_m^{ini}} \indi{\{w_p^{ini} < \wobjm \}} \mathbb{P}\left[u_{01}^p < \frac{w_p^{ini}}{\wobjm} \right]
+ \indi{\{w_p^{ini} \geq \wobjm \}} \mathbb{E}[N^{split}_p] \\
&= \frac{S_m^n}{S_m^n+E_{r,m}^n} N_{obj,m}^n 
+ \sum_{p = 1}^{N_m^{ini}} \frac{w_p^{ini}}{\wobjm}
= \frac{S_m^n}{S_m^n+E_{r,m}^n} N_{obj,m}^n + \frac{E_{r,m}^n}{\wobjm} = N_{obj,m}^n
\end{align}
since $\wobjm = \frac{S_m^n+E_{r,m}^n}{N^n_{obj,m}}$.

Moreover, following the specifications of the algorithm \ref{algo:wobj_homo} :
\begin{equation}
\mathbb{E}[N_{obj,m}^n] = \left\lfloor\dfrac{E_{r,m}^n + S_m^n}{E_{r,tot}^n}\right\rfloor N_{share}^n + \mathbb{P}[u_{01}^m < R_m] = \frac{E_{r,m}^n + S_m^n}{E_{r,tot}^n} N_{share}^n.  
\end{equation}

Also, let us denote by $\tilde{N}_{total}^n$ the total number of particles at the end of the algorithm \ref{algo:wobj_homo}; invoking the properties of the conditional expectation:
\begin{equation}
\begin{split}
\mathbb{E}\left[\tilde{N}_{total}^n\right] 
&= \sum_{N_{obj,m} = \indi{E_{r,m}^n >0} + \indi{S_m^n>0}} \indi{E_{r,m}^n >0} + \indi{S_m^n>0} +  \sum_{N_{obj,m} \neq \indi{E_{r,m}^n >0} + \indi{S_m^n>0}} \mathbb{E} \left[\tilde{N}_{m}^n \right] \\
& \le 
2 c_{\mathcal{M}} +  \sum_{N_{obj,m} \neq \indi{E_{r,m}^n >0} + \indi{S_m^n>0}}  \mathbb{E} \left[   \mathbb{E}[\tilde{N}_m^n|N_{obj,m}^n] \right] \\
& \le 2 c_{\mathcal{M}} +  \sum_{N_{obj,m} \neq \indi{E_{r,m}^n >0} + \indi{S_m^n>0}} \mathbb{E}[N_{obj,m}^n]  \\
& \le 2 c_{\mathcal{M}} + \sum_{N_{obj,m} \neq \indi{E_{r,m}^n >0} + \indi{S_m^n>0}}  \frac{E_{r,m}^n + S_m^n}{E_{r,tot}^n} N_{share}^n \le 2 c_{\mathcal{M}} + N_{share}^n = N_{total}^n.
\end{split}
\end{equation}

If $N_{m,obj}^n  = \left\lfloor\dfrac{E_{r,m}^n + S_m^n}{E_{r,tot}^n}\right\rfloor N_{total}^n + \indi{u_{01}^m < R_m  } > 0 \ \forall m \in \mathcal{M}$, then :
\begin{equation}
\mathbb{E}\left[\tilde{N}_{total}^n\right] 
= \sum_{m \in \mathcal{M}} \mathbb{E} \left[\tilde{N}_{m}^n \right], 
\end{equation}
and we conclude with the same argument as above.
\end{proof}

\begin{remark} Contrary to the algorithm that aims to obtain an uniform target number of particles by cell, the algorithm \ref{algo:wobj_homo} does not allow, 
	even if repeated several times, to have the certainty to obtain exactly  $N_{total}^n$ (or $N_{share}^n$) particles. Indeed, the algorithms in the section
	\ref{sec:algo_proposes} have the  advantage to converge exponentially towards $\nobj$. 
	However, algorithm \ref{algo:wobj_homo} allows  to control easily the number of the particles in the simulation. 
\end{remark}

Up to now, we presented three algorithms that improve the weight distribution, before the tracking phase.
The next part illustrates the effect of these algorithms on the weight distribution and its contribution for a radiative transfer example.

\section{Numerical results} \label{sec:res_num}

The \lcl \ algorithm influences the particle weight distribution before the tracking phase.
The subsection \ref{sec:vect_norme} will investigate the speed of convergence of the algorithm with or without conservative splitting. To this end, 
we analyze what is the result of the repetition of the algorithms on a vector of norm $1$ and study the distance to an uniform distribution 
and the number of weights at each iteration.
Then in the subsection \ref{sec:transfert_rad}, we apply the considered algorithms  on a radiative transfer case with the propagation of a Marshak wave, then we add a second wave of lower intensity, all with the  FCI2 code \cite{dattolo01}. In this case, we are interested in the speed of the waves and the variance of the radiative temperatures obtained.

\subsection{Illustration of the weight distribution \label{sec:vect_norme}}

We use the notations of the section \ref{sec:algo_proposes} and consider the evolution of the process $N^l$ related to a set of weights $\{w^l_i\}_{i \in 1, \dots, N^l}$.
We start by illustrating the evolution of this process when there is no source term (noted $S$, with $S=0$ therefore $N^{vol}=0$) then we illustrate the consequences of 
a non-null value of $N^{vol}_m$.

\subsubsection{Absence of source term $N^{vol} = 0$}

To illustrate the speed of convergence of the algorithms defined in the section \ref{sec:local_originale} and studied in the section \ref{sec:algo_proposes} we generate, 
according to a uniform law, weights between $0$ and $1$ then we perform a renormalization to obtain a uniform distribution of weights with norm equal to $1$.
We illustrate the results obtained by algorithms with conservative (denoted \textit{Split C}) or non-conservative splitting (denoted \textit{Split NC}) 
by presenting the evolution of the process $ (N^l)_{l \in \mathbb{N}} $ depending on the number $l$ of repetitions of the algorithm \ref{algo:meth_lcl} 
and the distance of the distribution $\{w^l_i \}_{i \in {1,\dots,N^l}} $ to a uniform distribution denoted $d^l$ defined simply by:
\begin{equation}
d^l = \sum_{i = 1}^{N^l} \left| w_i^l - \frac{1}{\nobj} \right|.
\end{equation}

\begin{figure}[H]
	\centering
	\begin{tabular}{cc}
		\includegraphics[scale=0.9]{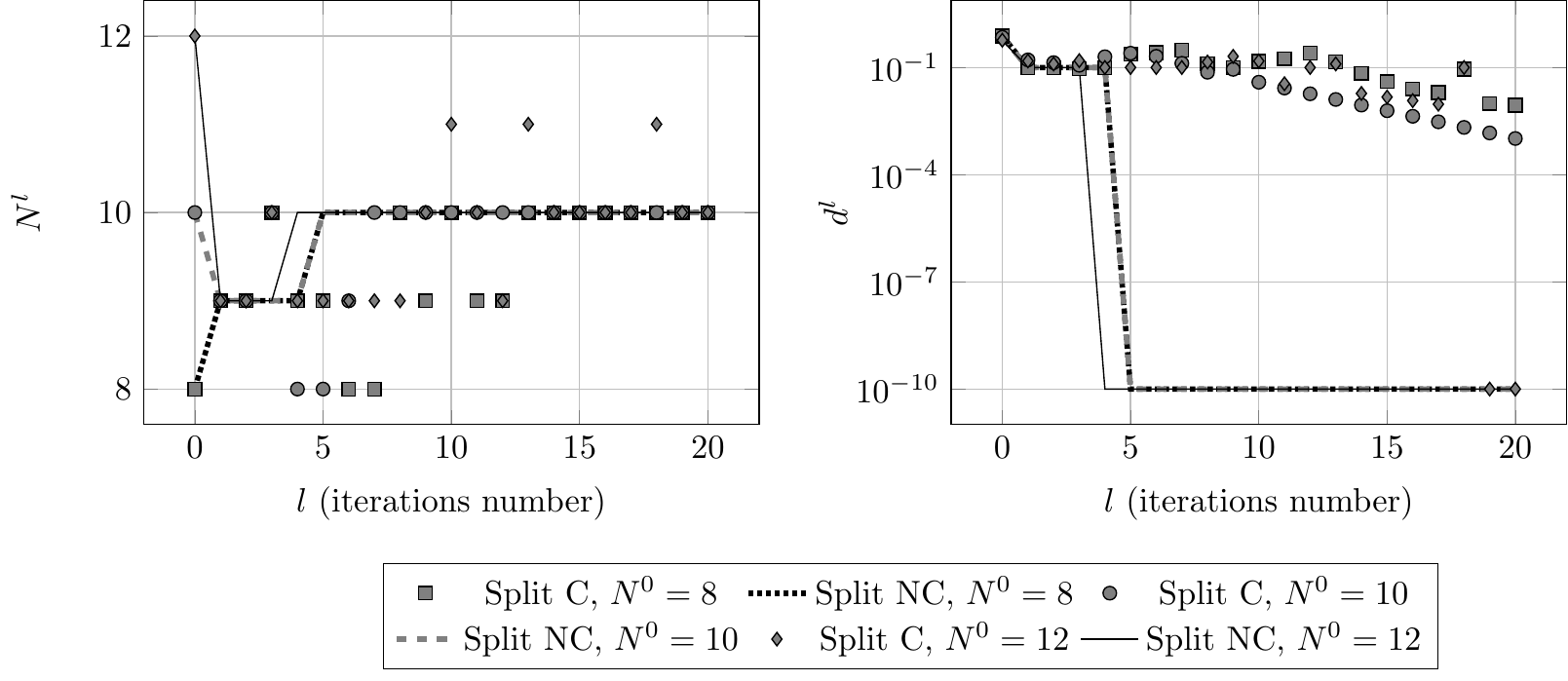}
	\end{tabular}
	\caption{Results obtained for $\nobj = 10$ with the algorithm \ref{algo:meth_lcl} and $S= 0$ on a population of weights of  norm $1$.
	\label{fig:source_0p0_nobj=10}}
\end{figure}

\begin{figure}[H]
    \centering
    \begin{tabular}{cc}
        \includegraphics[scale=0.9]{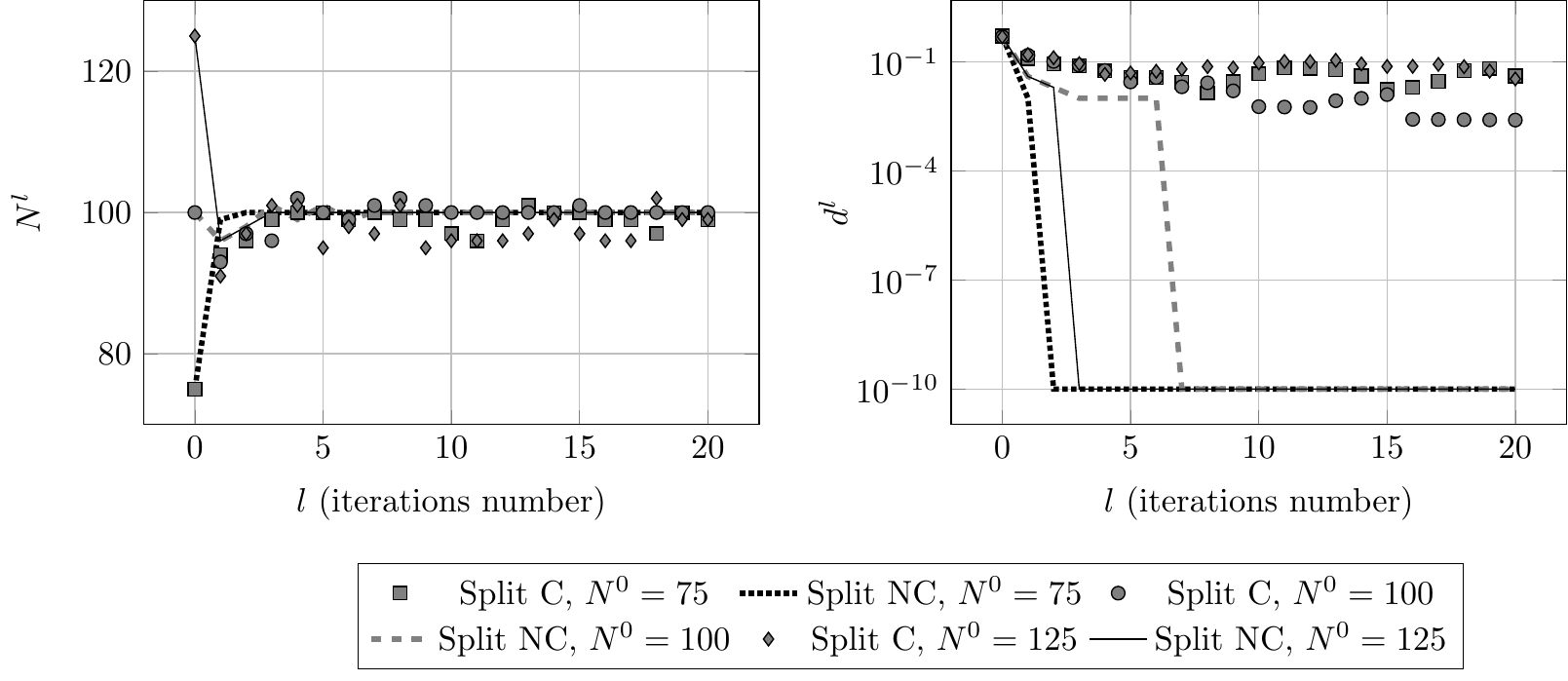}
    \end{tabular}
    \caption{Results obtained for $\nobj = 100$ with the algorithm  \ref{algo:meth_lcl} and $S= 0$  on a population of weights of  norm $1$. The difference between the conservative
    	and non-conservative versions is larger than in figure \ref{fig:source_0p0_nobj=10}.
    	\label{fig:source_0p0_nobj=100}}
\end{figure}

The figures \ref{fig:source_0p0_nobj=10} and \ref{fig:source_0p0_nobj=100} present the results obtained for different values of $\nobj$ and $N^0$.
The algorithm with \textit{Split NC} converges in a few iterations in both cases while the \textit{Split C} version is slower (even more when $\nobj$ increases).
Moreover, these curves illustrate that the state $N^l = \nobj$ is absorbing for the \textit{Split NC} algorithm while this is not the case for the \textit{Split C} 
counterpart who arrives there several times and leaves.

\subsubsection{Presence of a source term: $N^{vol}>0$}

We illustrate the algorithms in section \ref{sec:algo_proposes} with a source term in a manner analogous to the previous section with the exception that the 
generated weight renormalization is not at $1$ but at $1-S$ with $S \in [0, 1]$ representing the source term.
The algorithm \ref{algo:meth_lcl} is then used after which only the Russian-Roulette, Splitting and renormalization phases are repeated on the weights not coming from $S$ with the update of $\wobj$ as discussed in the section  \ref{sec:rq_terme_source}.

\begin{figure}[H]
	\centering
	\begin{tabular}{cc}
		\includegraphics[scale=0.8]{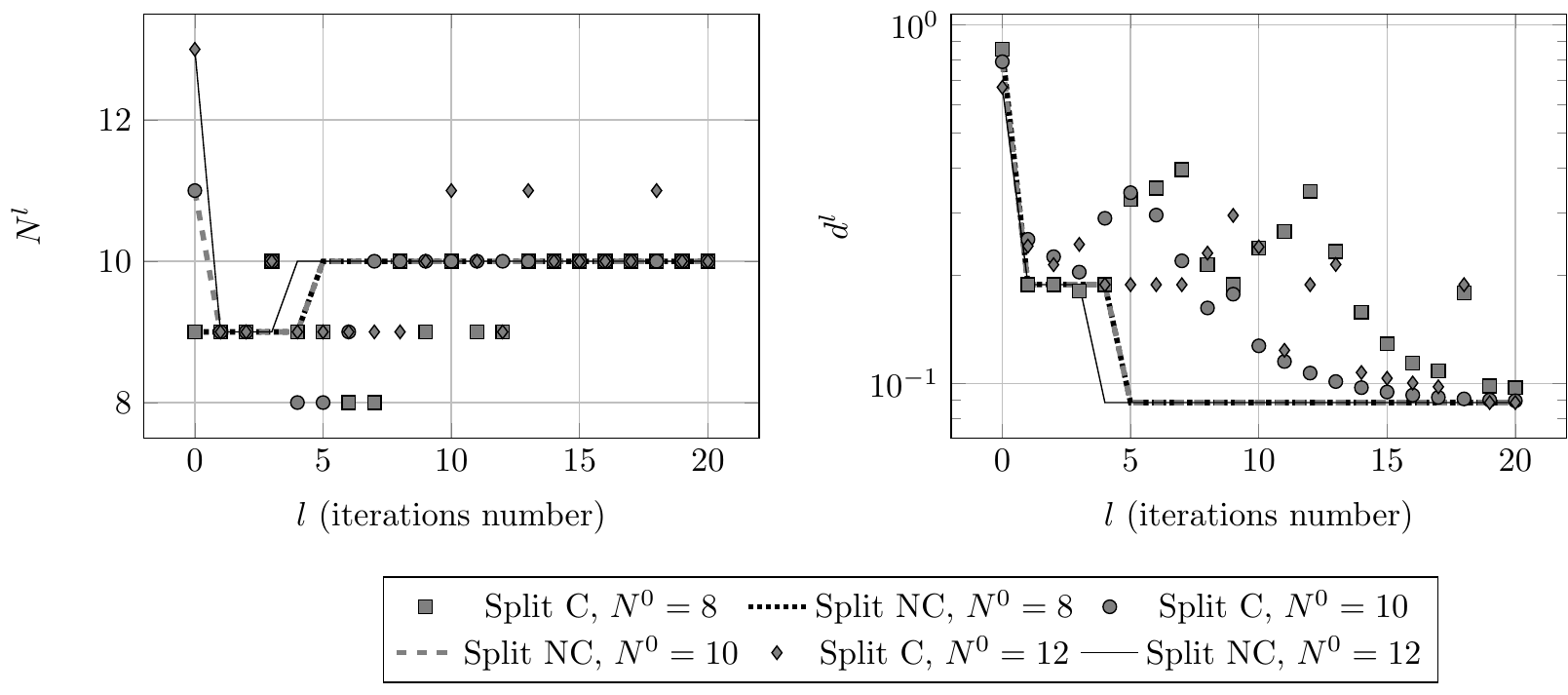}
	\end{tabular}
    \caption{Results obtained for $\nobj = 10$ with the algorithm  \ref{algo:meth_lcl} and $S = 0.0115$  on a population of weights of  norm $1$. \label{fig:source_1p15_nobj=10}}
\end{figure}

\begin{figure}[H]
	\centering
	\begin{tabular}{cc}
		\includegraphics[scale=0.8]{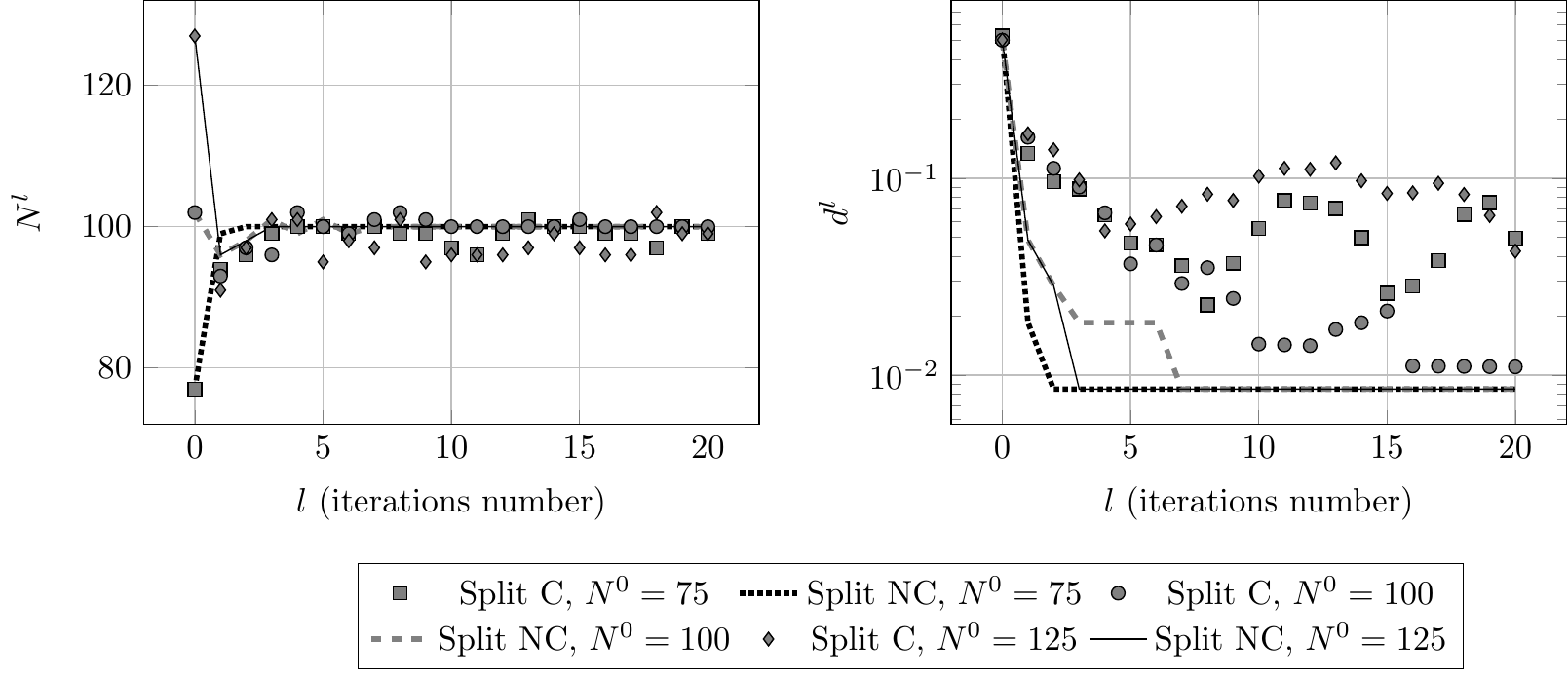}
	\end{tabular}
    \caption{Results obtained for $\nobj = 100$ with the algorithm  \ref{algo:meth_lcl} and $S = 0.0115$  on a population of weights of  norm $1$. \label{fig:source_1p15_nobj=100}}
\end{figure}

\begin{figure}[H]
	\centering
	\begin{tabular}{cc}
		\includegraphics[scale=0.8]{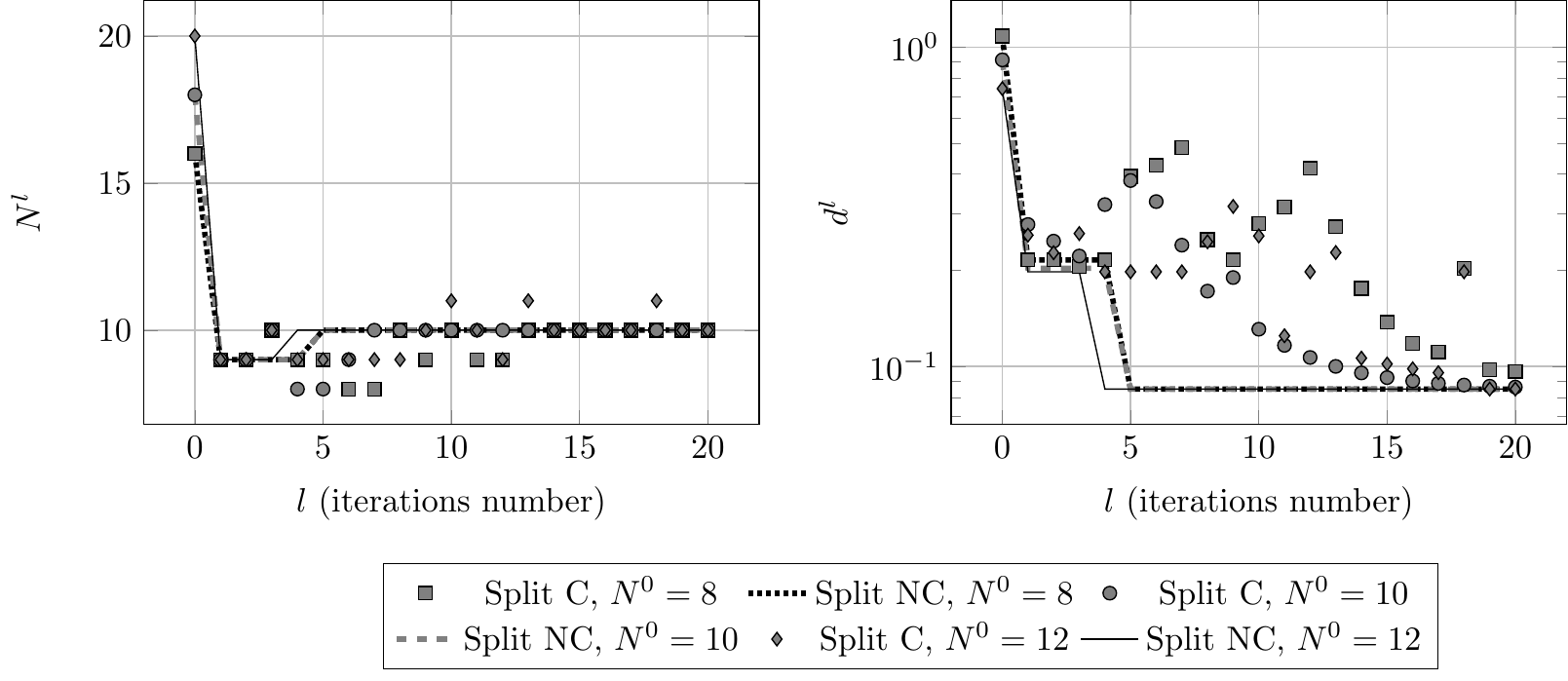}
	\end{tabular}
    \caption{Results obtained for $\nobj = 10$ with the algorithm  \ref{algo:meth_lcl} and $S = 0.715$  on a population of weights of  norm $1$. \label{fig:source_71p5_nobj=10}}
\end{figure}

\begin{figure}[H]
	\centering
	\begin{tabular}{cc}
		\includegraphics[scale=0.8]{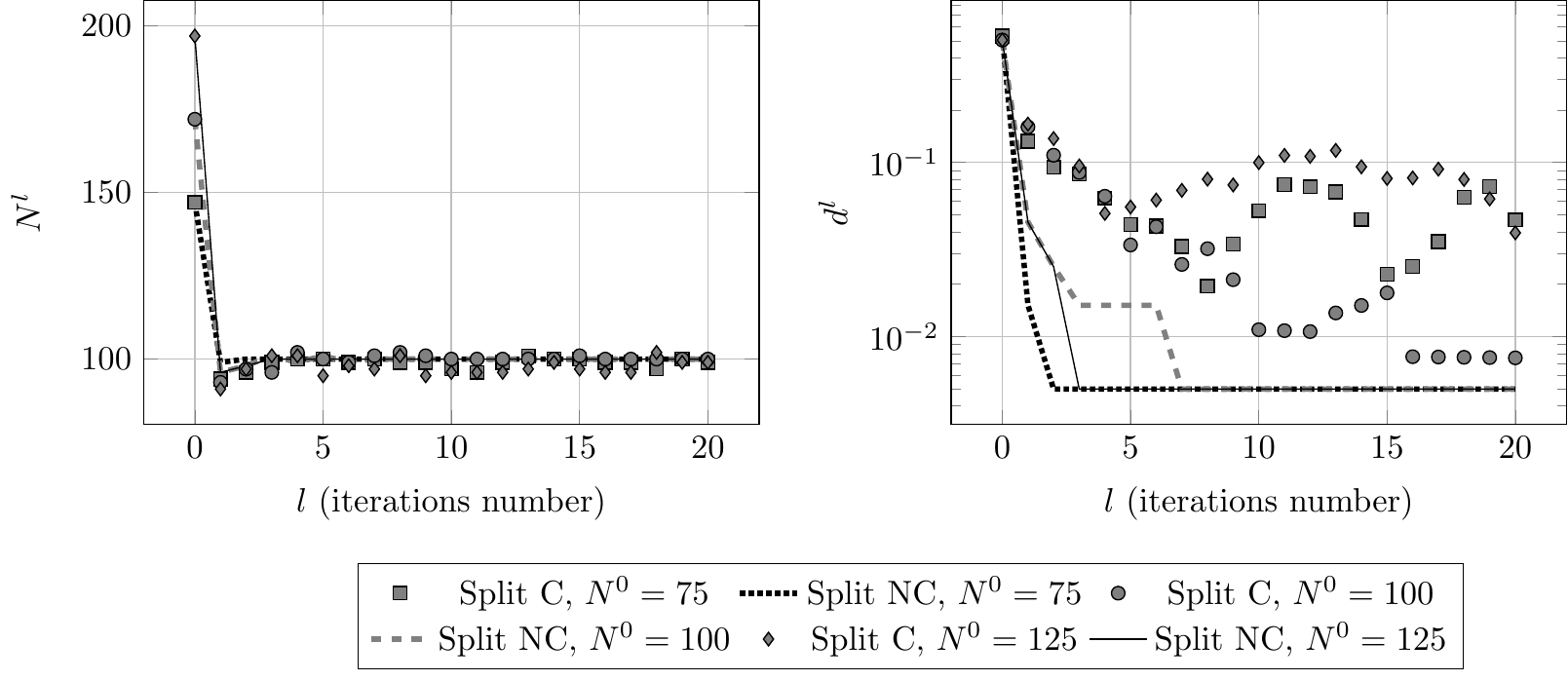}
	\end{tabular}
    \caption{Results obtained for $\nobj = 100$ with the algorithm  \ref{algo:meth_lcl} and $S = 0.715$  on a population of weights of norm $1$. \label{fig:source_71p5_nobj=100}}
\end{figure}

For a low value of $S$, the results obtained are shown in figures \ref{fig:source_1p15_nobj=10} and \ref{fig:source_1p15_nobj=100} and illustrate the influence of $ \nobj$ and $N_0$. 
Note that we have limited the minimum error to $10^{-10}$.
As in the previous case, the  \textit{Split NC} algorithm converges in a few iterations, unlike the \textit{Split C} algorithm.
However, the figure \ref{fig:source_1p15_nobj=10} shows that the limit distribution  is not necessary a uniform distribution: the presence of a source term, 
weak and much lower than $ \wobj $, limits the converges of $d^l$.
When the value of $\nobj$ increases (figure \ref{fig:source_1p15_nobj=100}) the distance to the uniform distribution decreases.
When the value of the source term is high, these findings remain unchanged as shown in figures \ref{fig:source_71p5_nobj=10} and \ref{fig:source_71p5_nobj=100}.

{\modif
\begin{remark}
For radiative transfer applications, the number of particles in figures~\ref{fig:source_0p0_nobj=10}, \ref{fig:source_1p15_nobj=10} and 
\ref{fig:source_71p5_nobj=10}
and  may appear to be small. But we think that it is useful to illustrate that the difference in convergence speed between conservative and non conservative splitting does not depend on the number 
of particles. In fact, when the number of cells is high (e.g. for three dimensional simulations) the number of particles per cell may be small.
\end{remark}
}

We have illustrated the results obtained by algorithms with conservative or non-conservative splitting in a simplified framework neglecting in particular the spatial and angular aspect.	The following part allows to apply the algorithms in a framework of the resolution of a transport equation via the equations of the radiative transfer.

\subsection{Radiative transfer: propagation of a Marshak-type wave \label{sec:transfert_rad}}

The algorithm \ref{algo:meth_lcl}  in its conservative or non-conservative splitting versions is compared with the proposed algorithm  \ref{algo:wobj_homo}.
The latter shows promising numerical results both in the context of the propagation of a Marshak-type wave and in the case of two waves having  intensities of different scale, 
in particular by reducing the calculation time and the overall variance of the system in the second case.

We apply the algorithms presented previously in the framework of the resolution of radiative transfer equations.
{\m We study the statistical noise for a fixed mesh and time step. We do not provide a global convergence study (in time, angle or space) as the calculation of the full, exact, solution is not the goal of this work; instead the Marshak wave is a well-know test case with documented reference solutions.}
The section \ref{sec:marshak_uo} is concerned  with the framework of the propagation of a Marshak-type wave.
Then the section \ref{sec:marshak_do} proposes, in order to challenge the influence of the spatial distribution proposed by the algorithm \ref{algo:wobj_homo}, a second wave of much weaker intensity than the first.

\subsubsection{Propagation of a Marshak wave \label{sec:marshak_uo}}

We test the different algorithms presented on the propagation of a Marshak-type wave in an opaque medium (see \cite{marshak1958,MCCLARREN20089711} for details)
using the FCI2 code described in \cite{dattolo01}.

We assume an ideal gas equation under the gray approximation.
The Monte Carlo method used here is based on the Fleck \& Cummings \cite{FLECK1971313} linearization mentioned in the section \ref{sec:math_frame}.
We use a model with two temperatures (radiative and matter): except  mention of the contrary, the term \textit{temperature} (noted $ T_{matter} $) will indicate the 
matter temperature since it is the one that appears in the emission term of the equation \eqref{eq:emission_volumique}.

\begin{figure}
	\centering
	\includegraphics[scale=1.0]{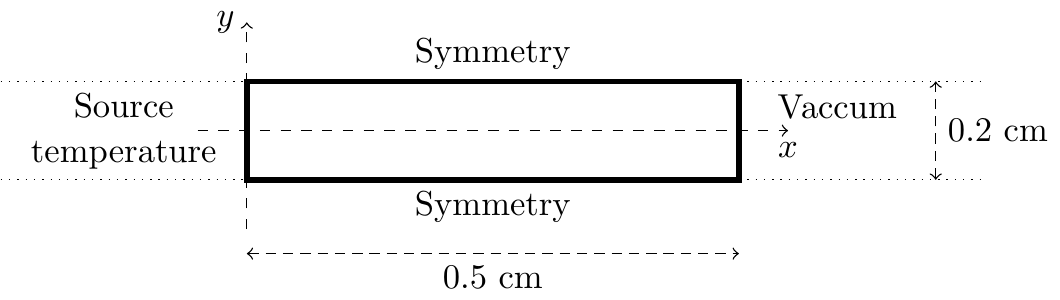}
	\caption{Diagram of the setting of the Marshak wave discretized into $50$ identical cells of size $0.01 \times 0.5$ cm. \label{fig:M1o_maillage}}
\end{figure}

This is a 1D benchmark (but treated in 2D with symmetry conditions on the top and bottom edges of the mesh (see figure \ref{fig:M1o_maillage}).
We then solve the equation \eqref{eq:transport_grise} in the section \ref{sec:math_frame} for
$j(t,x,\omega) =  k(t,x,\omega) = \rho \times d \times T_{matter}^{-3}(t,x)$.  The values and units used are specified in the  table \ref{tab:valeursMarshak}.

\begin{table}[h!]
	\centering
	\begin{tabular}{|c|c|c|c|c|}
		\hline
	    $I$ & $erg . cm^{-2} . s^{-1}$ & $a$ & $ 7.56 \times 10^{-15} \ erg . cm^{-3} . K^{-4} $ \\ \hline
	 	$dt$ & $4 \times 10^{-11} \ s$ & $d$ & $1.56 \times 10^{23} \  K^3 . g^{-1} . cm^2$ \\ \hline
	 	$\rho$ & $3 \  g . cm^{-3}$ & $c$ & $3 \times 10^{10} \ cm . s^{-1}$ \\ \hline
        $T_{matter}$ & $K$ & $T_{matter}(0,\cdot)$ & $11604 \ K$ \\ \hline
        $C_V$ & $8.6177 \times 10^{7} \ erg.g^{-1}.K^{-1}$ & $T_{matter}(\cdot, \text{left border})$ & $11604000 \ K$	 \\ \hline	 
	\end{tabular}
\caption{Values and units used in the numerical simulation of the propagation of a Marshak-type wave in an opaque medium. \label{tab:valeursMarshak}}
\end{table}

For each cell $ m $, at time $t^n$, the value of the volume emission is given by the formula \eqref{eq:emission_volumique}
and the surface emission term is given by the following formula:
\begin{equation}
\frac{acT_{matter}^4(\cdot, \text{left border}) \Delta x \Delta t }{4}.
\end{equation}

These two terms are positive: the \lcl \ method can therefore be applied.
$
\
$

We analyze the wave profile at $74$ns using a time step of $\Delta t = 4 \times 10^{-11} s$.
To do this, we perform $n=30$ realizations for $\nobj \in \{20, 200, 2000\}$ and we compare the statistical variance obtained per cell.
Note that in this case, the matter temperature is coupled to the radiative temperature.

\begin{figure}
    \includegraphics[scale=0.97]{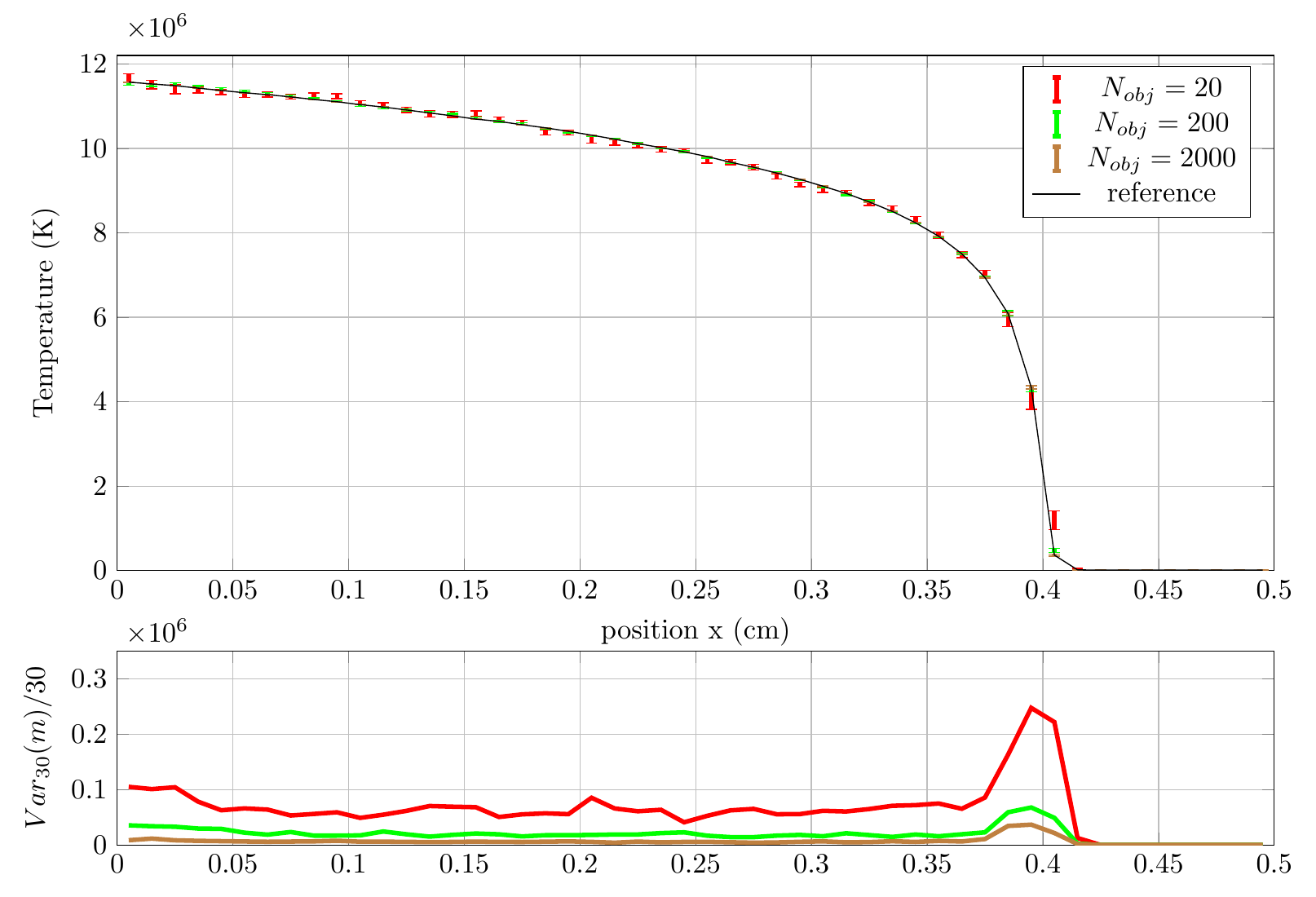}
	\caption{Results obtained with the algorithm \ref{algo:meth_lcl} and conservative splitting. Top: matter temperature. Bottom: variance per cell over $30$ realizations. The reference is the mean value of the temperature for the simulation with $\nobj = 2000$. \label{fig:Muo_tmCEStrue}}
\end{figure}

\begin{figure}
	\includegraphics[scale=0.97]{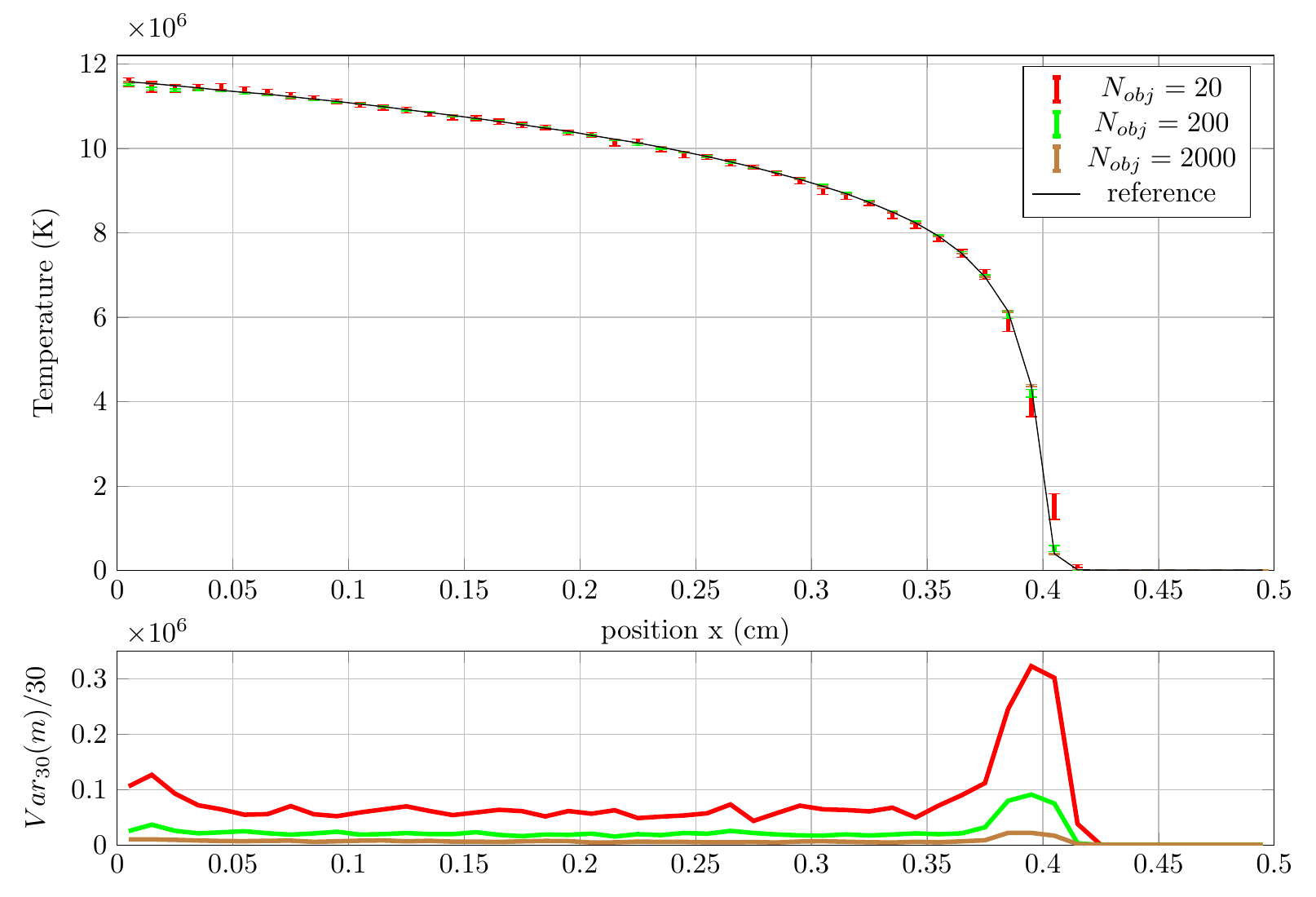}
	\caption{Results obtained with the algorithm \ref{algo:meth_lcl} and non conservative splitting. Top: matter temperature. Bottom: variance per cell over $30$ realizations. The reference is the mean value of the temperature for the simulation with $\nobj = 2000$.	
		\label{fig:Muo_tmCESfalse}}
\end{figure}

\begin{figure}
	\includegraphics[scale=0.97]{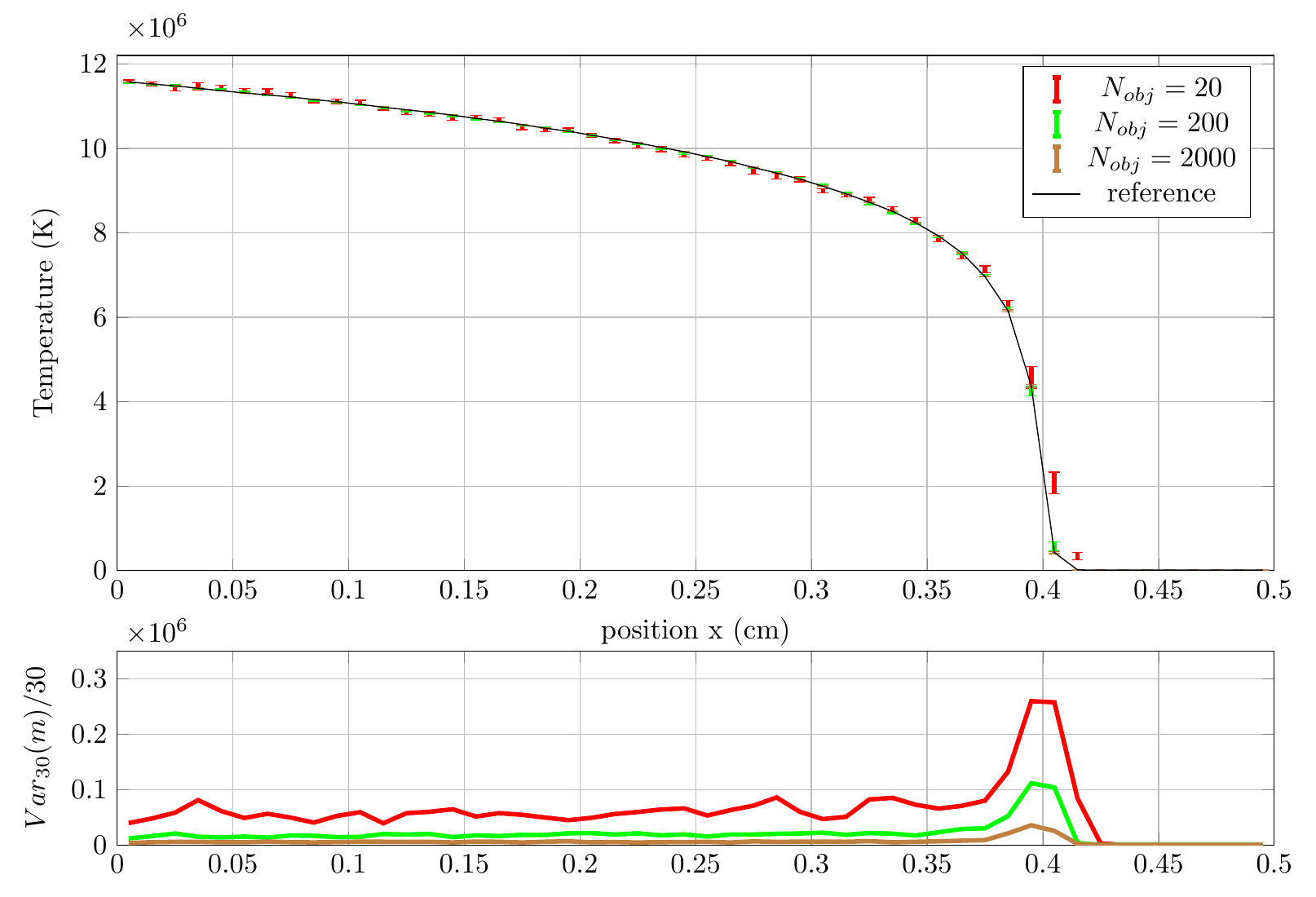}
	\caption{Results obtained with the algorithm \ref{algo:wobj_homo} and non conservative splitting. Top: matter temperature. Bottom: variance per cell over $30$ realizations. The reference is the mean value of the temperature for the simulation with $\nobj = 2000$.	
		\label{fig:Muo_tmWobjHomo}}
\end{figure}

\begin{figure}
	\centering
	\includegraphics[scale=0.97]{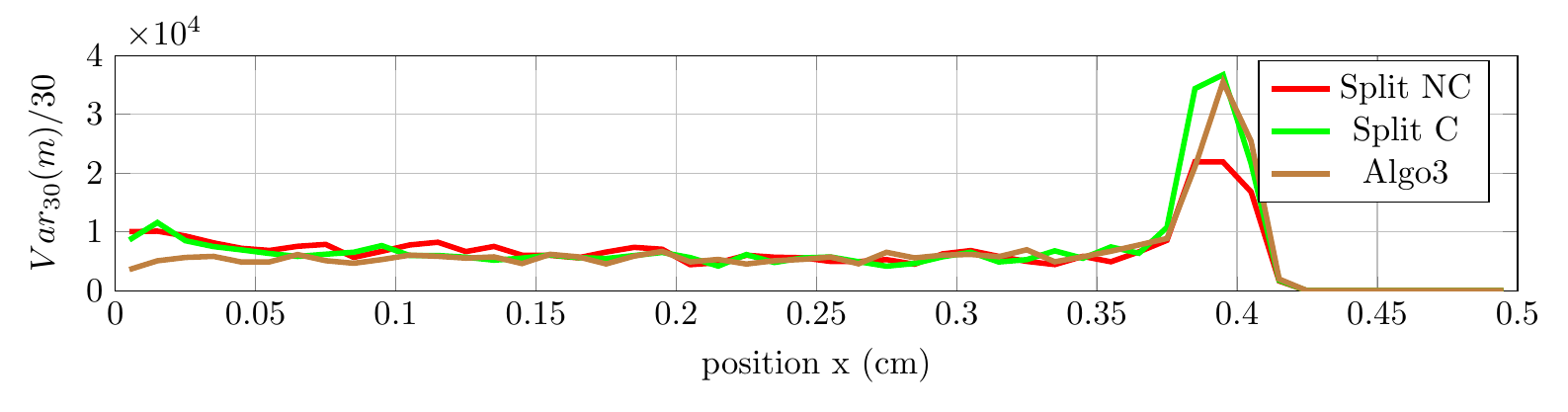}
	\caption{Comparison of the variance obtained with the three algorithms for $\nobj = 2000$. \label{fig:Muo_comparaison}}
\end{figure}

\begin{figure}
	\centering
	\includegraphics[scale=1.0]{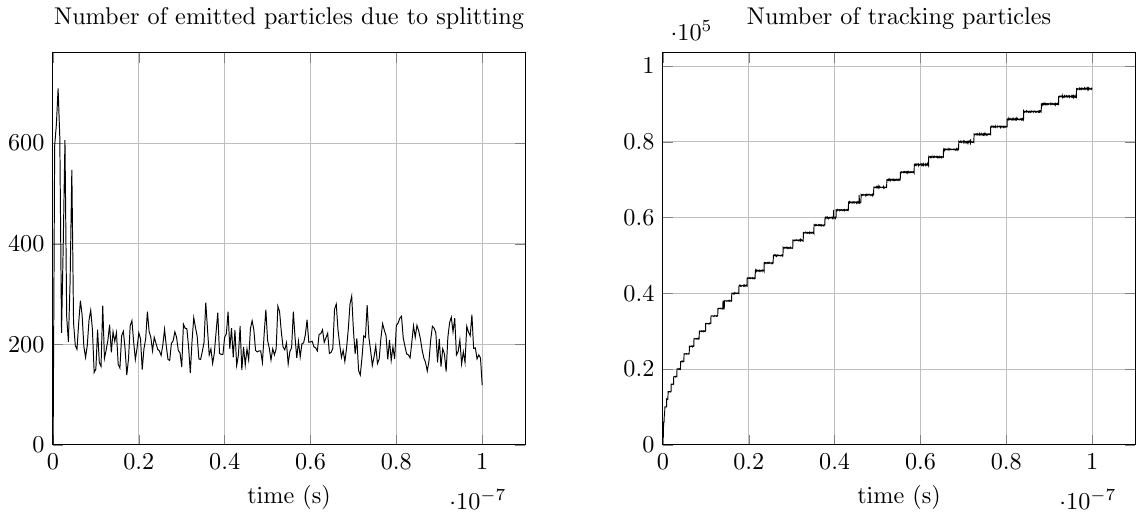}
	\caption{Number of particles {\md split} (left) and {\md number of individual tracks through the cell over the simulation} (right) over time in the case of the propagation of a Marshak type wave and $\nobj=200$. \label{fig:Muo_nbPartSplit_vs_Total}}
\end{figure}

\begin{figure}
	\centering
	\includegraphics[scale=1.25]{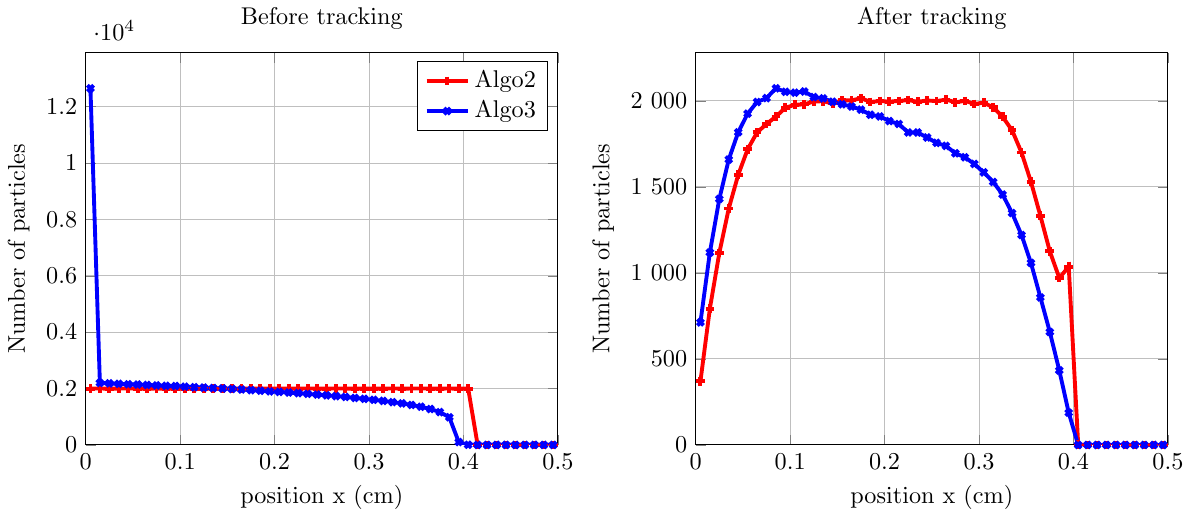}
	\caption{Distribution of particles on the cell before (left) and after (right) the last tracking phase in the case of the propagation of a Marshak wave.
		\label{fig:Marshak_density_track}}
\end{figure}

The variance per cell of the results obtained and the temperature profile with a confidence interval of 99\% are the object of the figures \ref{fig:Muo_tmCEStrue} and \ref{fig:Muo_tmCESfalse} when the objective number per cell is homogeneous (algorithm \ref{algo:meth_lcl}) and the splitting is conservative or not;  
the figure  \ref{fig:Muo_tmWobjHomo} presents the results of the application of the algorithm \ref{algo:wobj_homo}.
For the three methods, the cell at the foot of the wave has the greatest variance.
The figure \ref{fig:Muo_comparaison} allows to compare the size of the confidence interval per cell: it is interesting to note that non-conservative splitting does improve the variance at the bottom of the wave, even if there are very few splitted particles compared to the total number of particles tracked. The figure \ref{fig:Muo_nbPartSplit_vs_Total} shows the evolution over time of the number of splitted particles compared to the total number of particles of the simulation: only about 0.25\% of particles are affected.
The results obtained by the  algorithm \ref{algo:wobj_homo} are close to the case with conservative splitting at the wave foot rather than with non-conservative splitting;
nevertheless these results are improved on the left edge, the origin of the source temperature condition as shown in figure  \ref{fig:Muo_tmWobjHomo}.
This is explained by the density of particles: the figure \ref{fig:Marshak_density_track} shows the particle densities before and after the last tracking phase. The foot of the wave is depopulated with the algorithm \ref{algo:wobj_homo} while the limit condition is much better represented.

{\modiftwo To present the results we use a "figure of merit" metric.
Denote by $n$ the total number of realizations and $T_{matter}(u,m)$ 
the matter temperature on the cell 
$m \in \mathcal{M}$ at realisation $u$. The figure of merit, denoted FOM from now on, is  
 defined \cite{juzaitis80,mosher2009,wagner1998,Forster,mckinney2000} as the inverse of the product between the average CPU time  $T_{CPU}(n)$ and the square  relative error $RE^2(n)$ :
 \begin{equation}
 	FOM(n) = \frac{1}{RE^2(n) \times T_{CPU}(n)}.
 	\label{eq:FOMS}
 \end{equation}
 When the output of interest is a scalar, the relative error  $RE(n)$ is simply the standard deviation divided by the average value; however here the output (matter temperature) is a vector indexed over the cells $m \in \mathcal{M}$; accordingly, $RE^2(n)$ will be the norm (squared) of the vector of standard deviations divided by the norm (squared) of the average matter temperature vector. More precisely, to obtain the $RE^2$, we compute successively :
 
-  the average matter temperature vector  : 
 $T_{matter}(m) = \frac{1}{n} \sum_{u=1}^n T_{matter}(u,m)$, $m \in \mathcal{M}$;
 
- its squared norm :   $\| T_{matter}(\cdot) \|^2 = \sum_{m \in \mathcal{M}} T_{matter}(m)^2$;

- the unbiased empirical standard deviation matter temperature vector 
$ Std(m) = \sqrt{Var(m)}$ where 
$Var(m)= \frac{1}{n-1} \sum_{u=1}^n \left[ T_{matter}(u,m) - T_{matter}(m) \right]^2 $, $m \in \mathcal{M}$;

- its squared norm :   $\| Std(\cdot) \|^2 = \sum_{m \in \mathcal{M}} Std(m)^2 =  \sum_{m \in \mathcal{M}} Var(m)$

- and finally
\begin{equation}
RE^2(n) = \frac{\| Std(\cdot) \|^2}{\| T_{matter}(\cdot)\|^2} = \frac{\sum_{m \in \mathcal{M}} Var(m)}{\sum_{m \in \mathcal{M}} T_{matter}(m)^2}.
	\label{eq:REsq}
\end{equation}

}

\begin{figure}
	\includegraphics[scale=0.9]{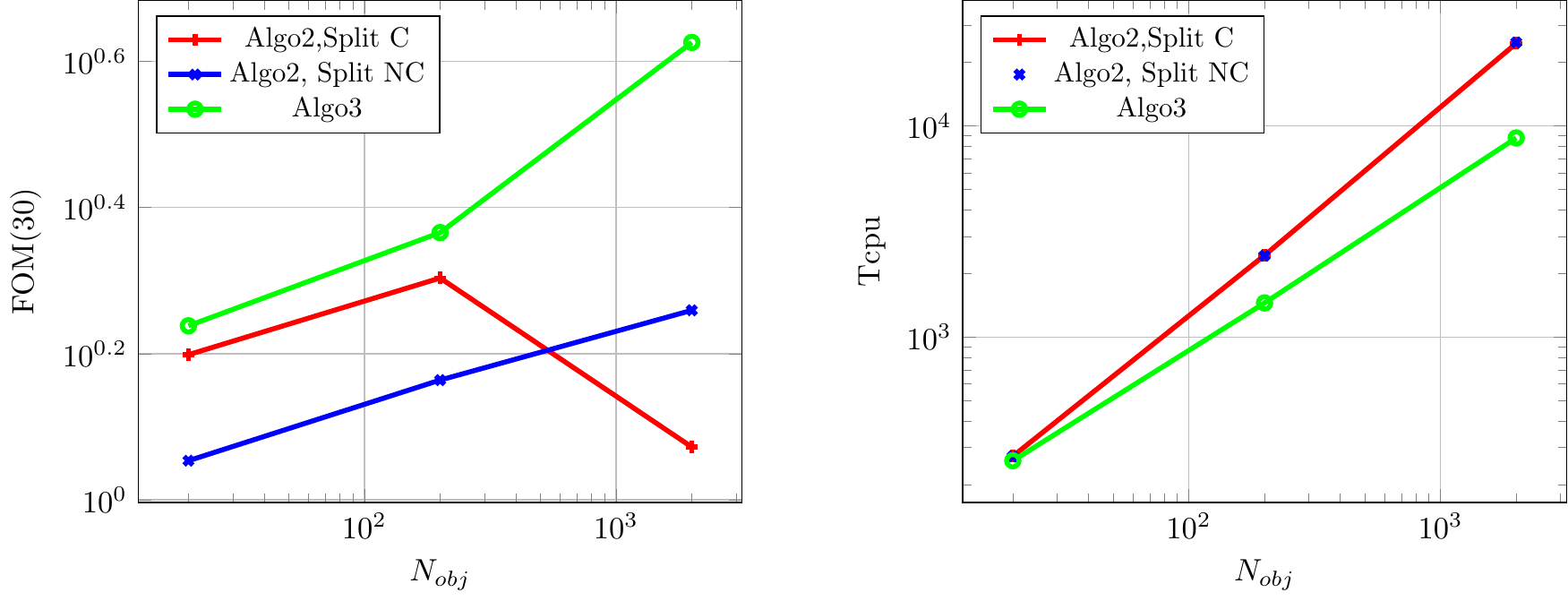}
	\caption{\modiftwo Metric $FOM(30)$ (left) and CPU time (right) for several values of $\nobj$ in the case of the propagation of a Marshak wave. 
	\label{fig:Muo_diagIC_fom}}
\end{figure}

The figure  \ref{fig:Muo_diagIC_fom} shows that the algorithm \ref{algo:wobj_homo} does not allow a significant gain in variance for the same number of particles followed but allows a significant gain in terms of computation time.
This case has the particularity of being placed at the diffusion limit: one explanation is that the particles placed at the foot of the wave have a higher CPU cost because they enter a very opaque medium (therefore undergo many events), unlike the left boundary of the domain (see figure \ref{fig:Marshak_density_track} for the distribution of particles in the cell before and after the last tracking phase).
The algorithm \ref{algo:meth_lcl} allows to have, on average, $\nobj$ particles per cell, both at the foot of the wave and at the left border.
It is interesting to note that having few particles at the foot of the wave, with the algorithm \ref{algo:wobj_homo}, does not prevent a variance comparable to the case with 
$\nobj$ particles and conservative splitting. In addition, the distribution after the tracking phase is similar.
Nevertheless, this depopulation of the foot of the wave could be a source of problems in the case of several distinct phenomena to be captured simultaneously: 
the following section provides an example in this context.

\subsubsection{Propagation of a primary and secondary wave \label{sec:marshak_do}}

The algorithm \ref{algo:wobj_homo} has shown its advantages within the framework of the propagation of a Marshak type wave: thanks to a spatial distribution of the particles, it allows, with comparable variance, to reduce the computation time.
However, an important question in the context of a multi-physics code, where the solution must be computed with enough precision over the entire domain,
is to make sure not to depopulate one area of the field in favor of another area and therefore 
to fail capturing the start of a phenomenon, however small it may be.
This is why we are considering the previous numerical application but adding a wave, of much lower intensity, on the right edge of the domain.
Moreover, on the right part of the domain, we go beyond the framework of the diffusion limit by modifying the emission and the absorption opacity by the values specified in the table 
\ref{tab:valeursMarshak_droite} in order to describe two phenomena which are different in both intensity and propagation speed. 
We observe the position of the two waves at $2$ns, this time using the radiative temperature, which better reflects the propagation of the low intensity wave.

\begin{remark}
It would have been possible to adapt the emission and the absorption opacity of the right part to remain in the diffusion limit. 
However, this would have led to expensive calculations for the right wave to propagate given its low intensity.
It is the same for the choice of the observed temperature: it is possible to consider the matter temperature, but in this case, it is necessary to increase the size of the domain to allow the second wave to propagate significantly without encountering the first wave.
This leads to a drastic increase of the final time of the simulation, the size of the cell and the number of particles followed thus increasing its cost.

\end{remark}

\begin{table}[h!]
	\centering
	\begin{tabular}{|c|c|c|c|c|}
		\hline
		$d_{right}$ & $1.56 \times 10^{13} \  K^3 . g^{-1} . cm^2$  & $T_{matter}(\cdot, \text{right border})$ & $116040 \ K$	 \\ \hline	 
	\end{tabular}
	\caption{Values used in the numerical simulation of a low intensity wave, to be compared with the table \ref{tab:valeursMarshak}. \label{tab:valeursMarshak_droite}}
\end{table}

\begin{figure}
	\includegraphics[scale=0.96]{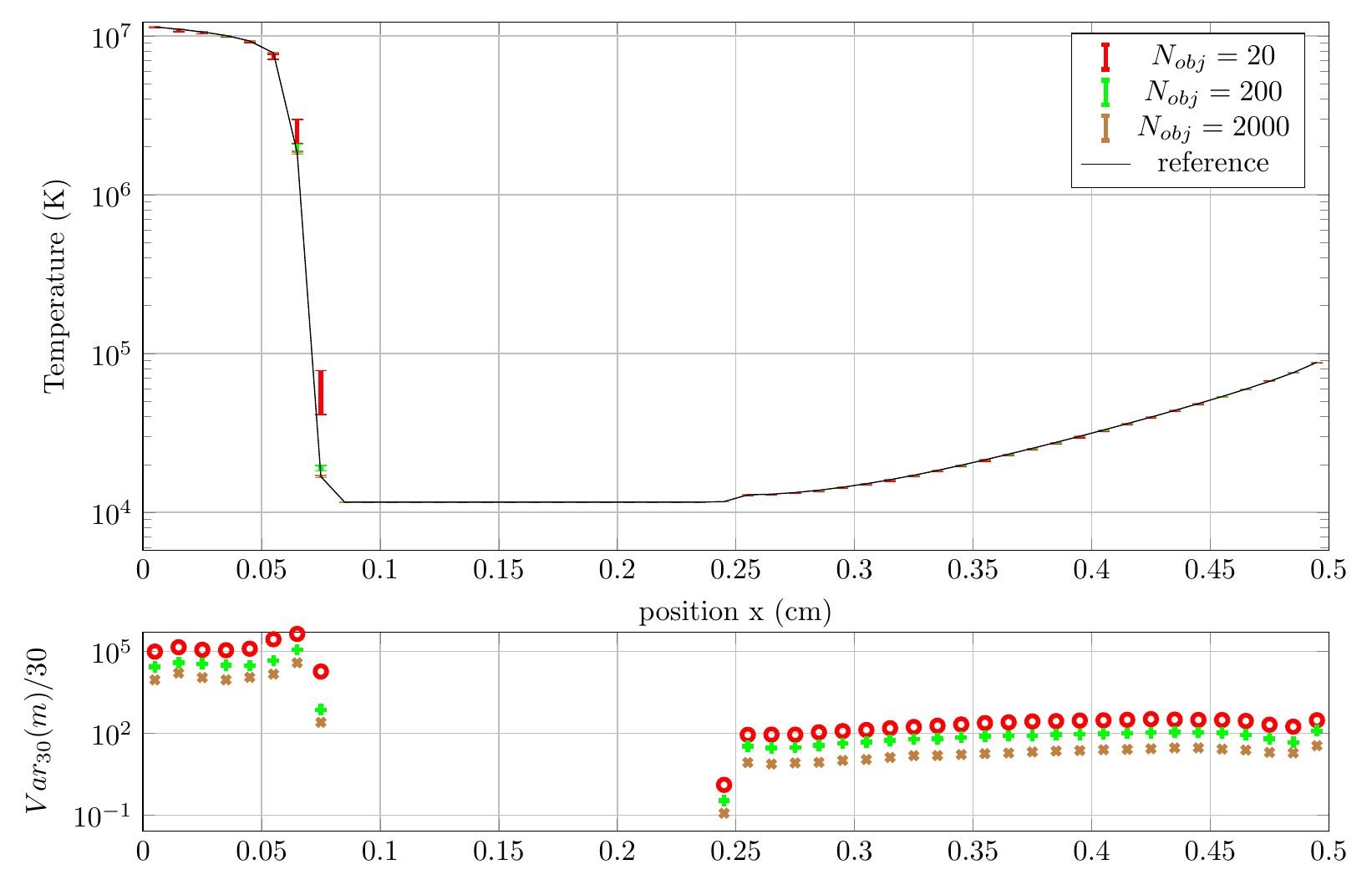}
	\caption{Results obtained with the algorithm \ref{algo:meth_lcl} and non-conservative splitting. Top: matter temperature. Bottom: variance per cell over $30$ realizations.  \label{fig:Mdo_tmCESfalse}}
\end{figure}

\begin{figure}
	\includegraphics[scale=0.96]{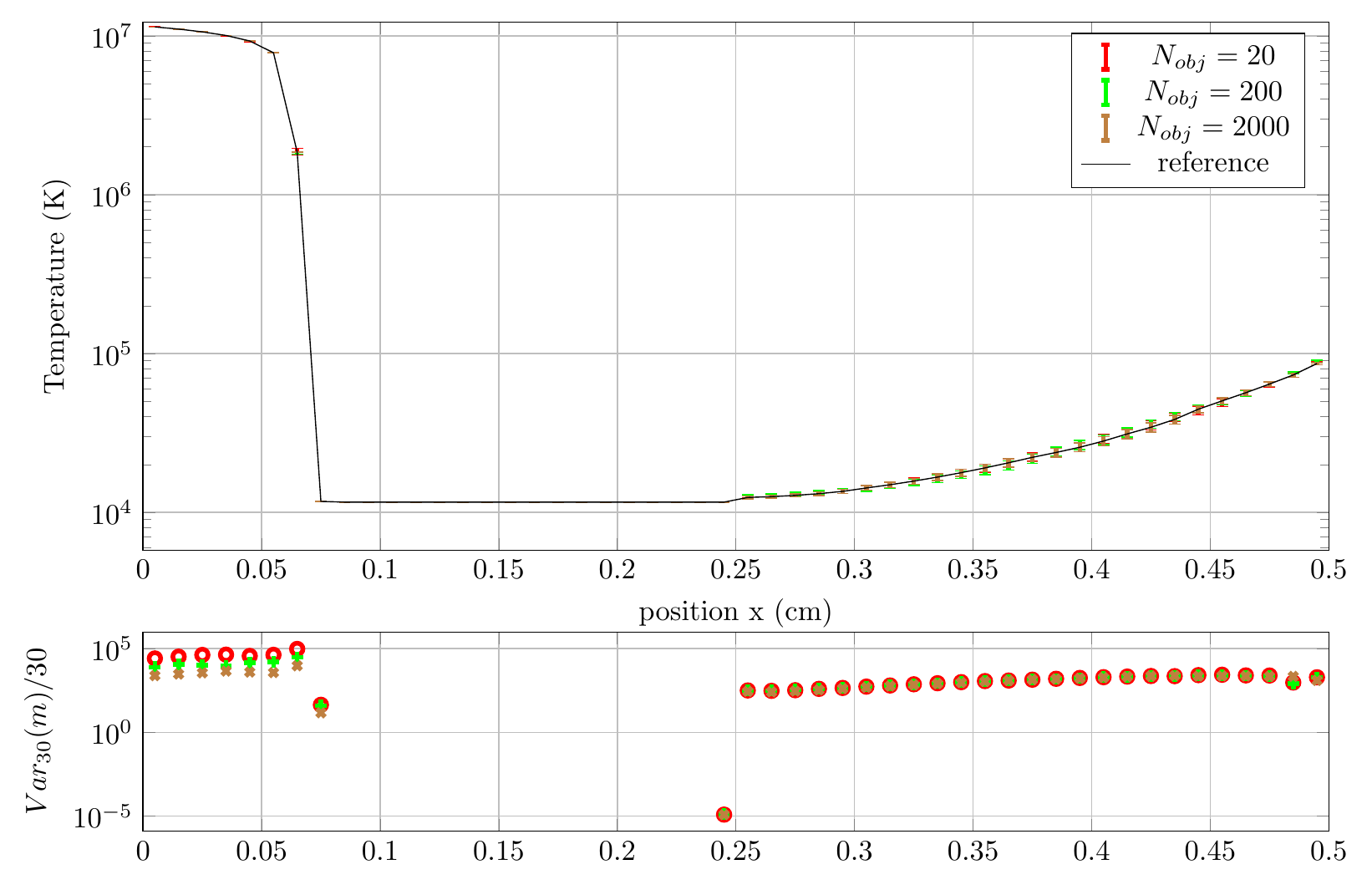}
	\caption{Results obtained with the algorithm \ref{algo:wobj_homo} and non-conservative splitting. Top: matter temperature. Bottom: variance per cell over $30$ realizations.  		
		\label{fig:Mdo_tmWobjHomo}}
\end{figure}

\begin{figure}
	\includegraphics[scale=1.250]{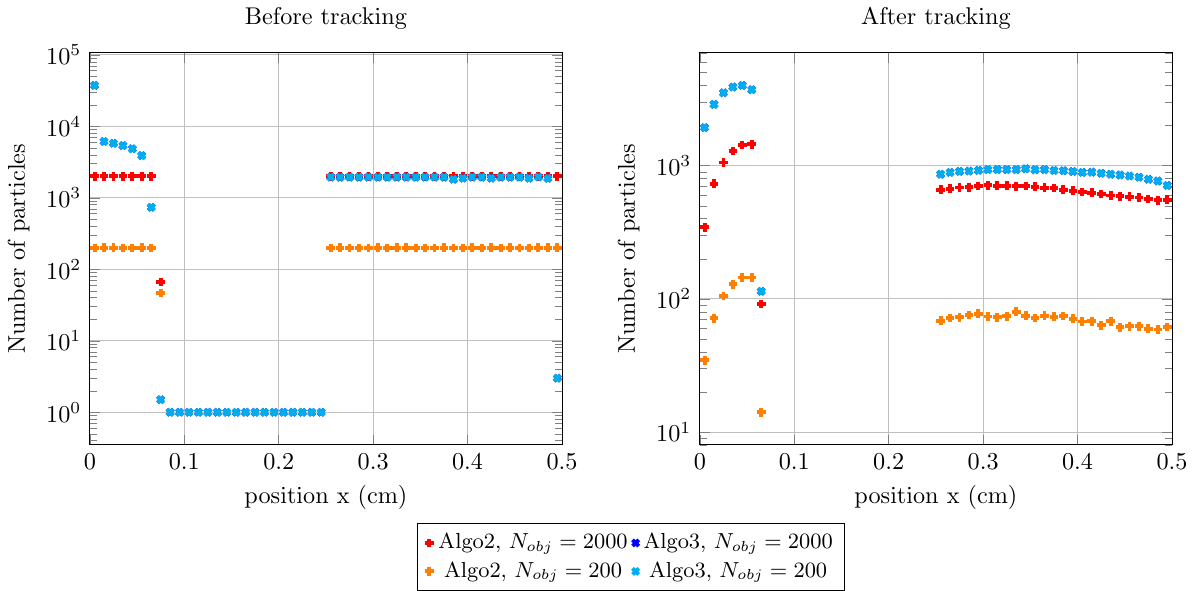}
		\caption{Distribution of the particles on the cell before (left) and after (right) the last tracking phase in the case of the propagation of a Marshak wave using the
	  algorithm \ref{algo:wobj_homo}  and different values of $\nobj$.\label{fig:M2o_density_track}}
\end{figure}

We present the same metrics as in the previous case, and the temperatures are presented in logarithmic scales so that the propagation of the weak wave can be represented.
Figure \ref{fig:Mdo_tmCESfalse} shows the results when $\nobj$ is homogeneous and the splitting is non-conservative; the figure  \ref{fig:Mdo_tmWobjHomo} 
treats the case when the algorithm \ref{algo:wobj_homo} is used.
The case with homogeneous $\nobj$ and conservative splitting is not presented in view of the results on the propagation of a single wave.
The propagation of the low intensity wave is reconstructed with both methods, even if with the algorithm \ref{algo:wobj_homo} few particles are used for this wave (figure 
\ref{fig:M2o_density_track}).
With the algorithm \ref{algo:wobj_homo}, increasing $\nobj$ does not change the variance on the low intensity wave, since the extra particles are allocated to the high intensity wave.
 
\begin{figure}
	\centering
	\includegraphics[scale=0.97]{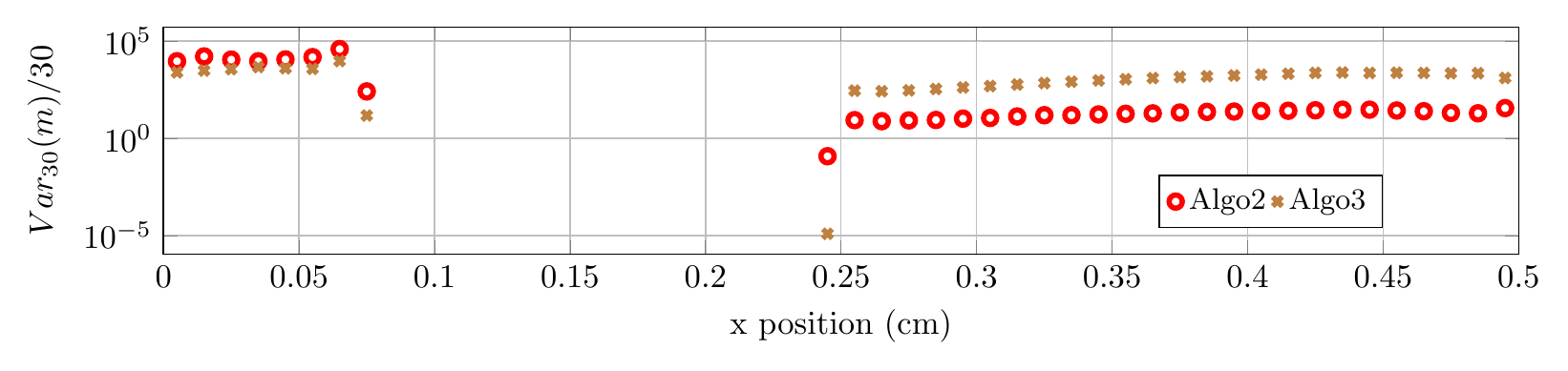}
	\caption{Comparison of the variance obtained with the two algorithms for $\nobj = 2000$. \label{fig:Mdo_comparaison}}
\end{figure}

Figures \ref{fig:Mdo_tmCESfalse} and \ref{fig:Mdo_tmWobjHomo} present the radiative temperatures and the associated confidence intervals and the figure \ref{fig:Mdo_comparaison} 
makes it possible to compare the size of the confidence intervals of the two methods for $\nobj=2000$. Unsurprisingly, the variance of the right wave is greater with the  algorithm 
\ref{algo:wobj_homo} than with the  algorithm \ref{algo:meth_lcl} (non-conservative splitting version).
The reverse occurs for the wave on the left. 
 
\begin{figure}
	\includegraphics[scale=0.85]{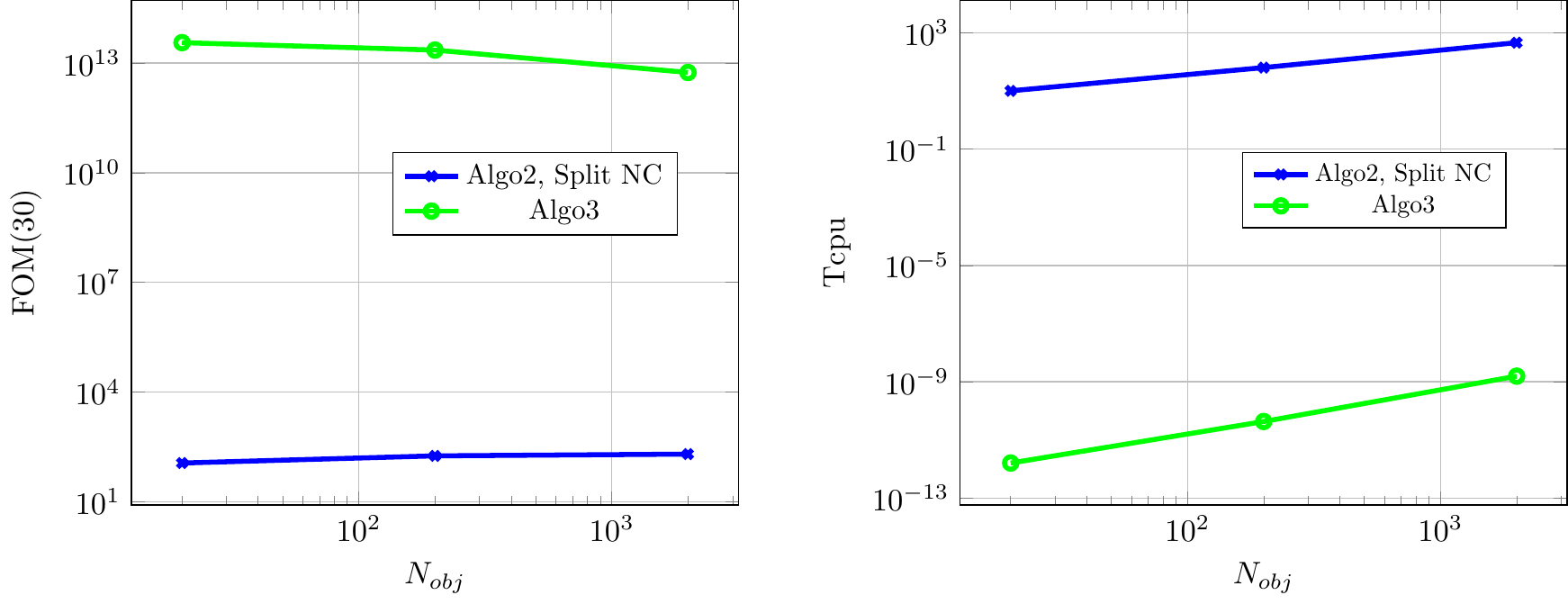}
	\caption{\modiftwo $FOM(30)$ metric (left) and CPU time (right) for several values of $\nobj$ for the case of the propagation of a Marshak wave with a second wave. 
		\label{fig:Mdo_diagIC_fom}}
\end{figure}

The figure \ref{fig:Mdo_diagIC_fom} {\modiftwo presents the FOM and the CPU time.}
The global variance is significantly reduced with the algorithm  \ref{algo:wobj_homo} which has a lower computation time than the algorithm \ref{algo:meth_lcl}. 
In this case, the use of a spatial distribution of the particles allows to gain both in global variance over the system and in computation time. 
This case underlines the main characteristic of the algorithm  \ref{algo:wobj_homo}  which aims to reduce the global variance on the domain.

This section allowed to challenge the algorithm  \ref{algo:meth_lcl} in a framework which seemed unfavorable. However, it made possible to underline the characteristics of this method through the reduction of the global variance in the system even if further practice is needed in order to better qualify these results.

\section{Conclusion}

The resolution of the transport equations by a Monte Carlo method allows a detailed level of modeling of the physical phenomena involved.
However, this stochastic method requires a precise sampling of the underlined Markov process in order to be able to treat the emission of the source term and the good 
representation of the solution obtained at the previous time step.
In this work, we have discussed the difficulties related to the variance reduction methods (cf., section \ref{sec:intro}) then we have presented the \lcl \  method which simultaneously performs the {\m control of existing population} phase and the source emission phase (cf., section \ref{sec:local_originale}).

First, based on the \lcl \ method presented in  \cite{laguzet2020}, we compared two types of splitting: first, one that conserves the  energy and secondly a
 non-conservative procedure (cf., section \ref{sec:algo_proposes}) (the total energy remains conserved in both cases thanks to a correction step).
These two algorithms, repeated several times on the same population of weights, converge towards the same limit distribution (cf., sections \ref{sec:algo_proposes} and 
\ref{sec:rq_terme_source}), but the algorithm with a non-conservative splitting has a faster convergence, as illustrated in section \ref{sec:vect_norme}.
Exploiting the properties of the \lcl \ method, we propose a variant in the section \ref{sec:lcl_wobj} which enables to obtain homogeneous weights over the entire cell considered. 
The algorithm \ref{algo:wobj_homo} contrasts with the initial method aiming to guarantee a homogeneous number of particles.

These algorithms are applied within the framework of a radiative transfer using the propagation of a wave in an opaque medium, corresponding to the diffusion limit of the system. 
For the three methods, the variance is largest at the foot of the wave. Non-conservative splitting, although involving a low proportion of particles in the system, allows a significant reduction in the variance at the foot of the wave. The newly proposed method allows to obtain a similar level of variance, but significantly reduces the calculation time and thus improves the merit index.

The algorithm  \ref{algo:wobj_homo}  is then tested on a situation which seems unfavorable: two physical phenomena distinct in intensity and speed are present on the considered domain.
In this context, the method with homogeneous objective number and non-conservative splitting is compared with the method aiming to obtain homogeneous objective weights. In this case, the variance on the main wave is improved to the detriment of the secondary wave. Nevertheless, at the global level of the simulation, the variance is improved, as well as the computation time.

$
\
$

The reduction of variance is a problem depending on the case studied.
It seems interesting to continue exploring this new method, in order to better characterize the situations where it is effective and the ones where it does not allow to obtain a satisfactory result on the whole domain (for instance, an area with a very important phenomenon but where energy is low). Finally, it is important not to fall back into too complicated  setups and parameter specifications: although the idea of defining an area where the weight should be homogeneous (for example distinguishing between two waves) is generally beneficial, it may also fall into the trap of the definition of a coefficient of importance; this is detrimental because the user is forced to employ the result of the calculation for the parameterization of the method.
Finally, in situations where zone depopulation can be too problematic, it would be interesting to combine the principle of the two methods: a minimum number of particles per cell, and distribute the rest in such a way as to homogenize the objective weights as much as possible.

$
\
$

\textit{This research did not receive any specific grant from funding agencies in the public, commercial, or not-for-profit sectors.}

\newpage

{\m
\appendix

\section{Proof of the proposition \ref{prop:cv_conservative} on the conservative splitting method}
\label{sec:appendix_proof_nc}

\textit{We provide in this appendix 
	a complete version of the proof sketched in main text; this version
 addresses rigorously all technical details required for the proof of
	 the proposition \ref{prop:cv_conservative} page \pageref{prop:cv_conservative}.
We recall that this proposition concerns one cell and a fixed time. In this appendix, we omit the index $m$ of a cell, and denote as before by $l$ 
the iteration counter.
}

$
\
$

Recall that $X_l$ is the state of the algorithm after $l$ iterations of (RR + S + the renormalization) steps. The $X_l$ have $N^l$ components noted $(X_l^1, \dots, X_l^{N^l}$).
Note that, in full generality, $X_l$ can be described as a set of (unknown number) positive weights that sum up to $E_{r,m}^l = 1$, i.e. a member of $\ell^1$, the ensemble of (absolutely) summable sequences. 
Denote also by $\target$ the uniform distribution with $\nobj$ particles all of same weight $\wobj$: 
$\target= (\wobj,...,\wobj)\in \R^{\nobj}$. 

$
\
$

Note for now that, since  $N_0$ is finite, $N_l$ will remain finite for any $l$; moreover, since for any $l$ the alternatives are in finite number (basically the outcomes of some binomial variables) then, given $X_0$, the set of all possible states taken by the algorithm is countable. 

As in the non-conservative case, the number of particles at the next step is depending of a sum of Bernoulli variables (noted $\textit{Be}$):
\begin{equation}
N^{l+1} =\sum_{i=1}^{N^l}  \left\lfloor\dfrac{X_l^i}{\wobj}\right\rfloor 
+ \textit{Be}\left(p_i \right).
\label{eq:Nnp1_conservatif} 
\end{equation} 
with $p_i = \frac{X_l^i}{\wobj} - \left\lfloor\dfrac{X_l^i}{\wobj}\right\rfloor$.

\noindent As $E_{r,m}^l =1$, $ \frac{1}{\nobj} \sum_{i=1}^{N^l} 
p_i \leq \sum_{i=1}^{N^l} X_l^i - \frac{1}{\nobj} \left\lfloor\dfrac{X_l^i}{\wobj}\right\rfloor\leq    \sum_{i=1}^{N^l} X_l^j \leq 1$ so  $ \sum_{i=1}^{N^l} 
p_i \leq \nobj \ \forall l \geq 1 $.

Therefore the algorithm generates a discrete time countable state Markov chain $(X_l)_{l\ge 0}$. We recall that such a stochastic process has the strong Markov property (we will use  this property later). Finally, note that, should the algorithm arrive at the ``non-void correction'' step, then the next distribution will consist of one particle with mass $E_{r,m}^l$ and the next one will be exactly $\target$. Therefore, without loss of generality we will consider only Markov chain scenarios that do not involve the non-void correction step.

$
\
$

With these preliminaries we can prove the proposition \ref{prop:cv_conservative}.
The proof requires many technical arguments and will be split into several parts. 

$
\
$

$\star$ The first step of the proof is to show that the Markov chain will arrive, after a finite time $\tau$, at a total number of particles $N_\tau \le \nobj$.
\begin{lemma}
	Under the hypothesis of Proposition \ref{prop:cv_conservative} there exists $c_1>0$ such that
	\begin{equation}
		\forall L \ge 0: \proba[ N_0 > \nobj, N_1 > \nobj,..., N_L > \nobj] \le e^{c_1 (L-1)}.
		\label{eq:returnnobj}
	\end{equation}
	\label{lemma:returntonobj}
	In particular, let us define a stopping time $\tau$ equal to the first $j\ge 1$ such that $N_j \le \nobj$. Then $\tau$ is finite (with probability one).
\end{lemma}
\begin{proof}
	The number of particles at step $l+1$ depends on the outcomes of a set of Bernoulli variables
	$Be(p_i)$ (equation \eqref{eq:Nnp1_conservatif}): in general their number can be arbitrary large, but the parameters $p_i$ (the means) are such that
	$p=\sum_i p_i$ is 
	an integer  
	less than $\nobj$.
	 For instance, if all but one particles are of weight less than $\wobj$ and the remaining particle of weight in $[2 \wobj,3\wobj[$ then
	$p$ will be $\nobj -2$. 
	
	When the random variables $Be(p_i)$ are such that
	$\sum_i Be(p_i)$ matches exactly its mean $p$ the next iteration $l+1$ will have exactly $\nobj$ particles.
	
	Let us recall the following inequality (known as the Chernoff bound \cite{Chernoff52,uspensky1937introduction,hagerup_guided_1990} for the Poisson binomial distribution $\sum_i Be(p_i)$) :
	\begin{equation}
		\proba\left[ \sum_i Be(p_i)  \ge p+1\right] \le e^{-1/(2+p)}.
	\end{equation}
	Since $p$ is bounded from below  
	(we excluded the situation $p=0$, see the discussion above) we obtain that there exists some constant $c_2 <1$ depending only on $\nobj$ such that for any $l\ge 1$: $\proba[ N_{l+1}  > \nobj ] \le c_2$ and also
	$\proba[ N_{l+1}  > \nobj | N_{l}  > \nobj] \le c_2$. 
	
	Then, for any $L\ge 1$:
	\begin{multline}
		\proba[ N_0 > \nobj, N_1 > \nobj,\dots, N_L > \nobj] = 
		\proba[N_L > \nobj| N_0 > \nobj, N_1 > \nobj, \dots, N_{L-1} > \nobj] \\
		\cdot
		\proba[N_{L-1} > \nobj| N_0 > \nobj, N_1 > \nobj, \dots, N_{L-2} > \nobj]  \cdot \dots
		\cdot \proba[N_1 > \nobj| N_0 > \nobj] \\ 
		= 
		\proba[N_L > \nobj| N_{L-1} > \nobj] \cdot
		\proba[N_{L-1} > \nobj| N_{L-2} > \nobj]  \cdot \dots
		\cdot \proba[N_1 > \nobj| N_0 > \nobj] \cdot \proba[N_0 > \nobj] \\
		 \le c_2^{L-1},
	\end{multline}
	where we used the Markov property to pass from  
	 $\proba[N_L > \nobj| N_0 > \nobj, N_1 > \nobj, \dots, N_{L-1} > \nobj]$
	to 
	$\proba[N_L > \nobj|N_{L-1} > \nobj]$.
	This gives \eqref{eq:returnnobj}  for  $c_1= -log(c_2) >0$.
	
	In particular the tail probability $\proba[\tau \ge L-1]$ decreases exponentially with $L$ and thus $\tau$ is finite. 
\end{proof}

$\star$ We invoke now the strong Markov property of the countable space discrete time Markov chain $(X_l)_{l \ge 0}$ to conclude that $(Y_l)_{l\ge 0} = (X_{\tau + l})_{l\ge 0}$ is also a Markov chain and moreover $Y_0$ has at most $\nobj$ non-null  weights (with probability one). Note that the conclusion 
in equation \eqref{eq:convergence_conservative} is true if and only if it is also true for the Markov chain $(Y_l)_{l\ge 0}$ because the limit is independent of what happens before the finite  stopping time $\tau$. Thus, without loss of generality and in order to keep the same  notations, we consider from now on that $X_0$ has at most $\nobj$ non-null weights.
Note that this does not necessarily imply the same for $X_1$ or following values of the Markov chain.

$
\
$

$\star$ We will need another result which is describing the stability of the target $\target$ with respect to the algorithm.
Consider 
$$\targete = \left \{ (w_1,...,w_{\nobj}) \in \R^{\nobj} |
\sum_{k=1}^{\nobj} w_k= \nobj \wobj,
\sum_{k=1}^{\nobj} |w_k-\wobj| \le \epsilon \wobj
\right\},$$
and denote $\tau_\epsilon$ the first iteration step $l$ when $X_l$ enters $\targete$.

\begin{lemma}
	There exists some constant $c_3 >0$ only depending on $\nobj$ such that
	\begin{equation}
		\proba[\lim_{l \to \infty} X_l = \target | \tau_\epsilon < \infty ] \ge 1-c_3 \epsilon.
	\end{equation}
	\label{lemma:stability}
\end{lemma}

\begin{proof}
	First note that we only need to prove the assertion for small enough $\epsilon$ because for large $\epsilon$ we have  $1-c_3 \epsilon \le 0$ and in this case the conclusion holds trivially because a probability is always a positive number. We can assume $\epsilon < 1$.
	Let us write 
	\begin{equation}
		\proba[\lim_{l \to \infty} X_l = \target | \tau_\epsilon < \infty ] 
		= \sum_{k \ge 1} \proba[\lim_{l \to \infty} X_l = \target | \tau_\epsilon =k ] 
		\proba[\tau_\epsilon =k | \tau_\epsilon < \infty ].
	\end{equation}
	Since $\sum_{k \ge 0}	\proba[\tau_\epsilon =k | \tau_\epsilon < \infty ]=1$ it is enough to find a constant $c_3$ such that for any $k$:   
	$$\proba[\lim_{l \to \infty} X_l = \target | X_k \in  \targete ]
	\ge 1-c_3 \epsilon,$$ where we used the Markov property of $(X_l)_{l\ge 0}$.

	As  $\epsilon < 1$, this means that any $X_k^j$ can be written as $X_k^j=\wobj(1+\epsilon_j)$ with $|\epsilon_j| \le \sum_j |\epsilon_j|=\epsilon <1$; moreover $\sum_j \epsilon_j=0$ and $\sum_j |\epsilon_j|=\epsilon$ imply 
	$\sum_{\{j; \epsilon_j \le 0 \}} |\epsilon_j|=\sum_{\{j; \epsilon_j \ge 0 \}} |\epsilon_j|=\epsilon/2$.
	
	When $\epsilon_j<0$  the probability for each  each $X_k^j$ to survive after RR is $1-|\epsilon_j|$ and when  $\epsilon_j>0$ the probability to survive as is (without any split) is 
	$1-\epsilon_j$. So the probability for each particle, being it subject to RR or S to survive and not undergo any split is at least $\prod_{j=1}^{\nobj} (1-|\epsilon_j|) \ge 1- \sum_j |\epsilon_j| = 1-\epsilon $. 
	Thus, with probability $1-\epsilon$ the distribution $X_{k+1}$ will have exactly $\nobj$ non null weights. 
	In this case, the renormalization factor (which will be less than unity because particles surviving RR increase the total mass) is exactly $\nobj/(\nobj+\epsilon/2)$.
	
	\noindent After the renormalization phase :
	\begin{itemize}
	\item the particles having $\epsilon_j \le 0$
	have survived after a RR phase, will be all equal to
	$\nobj/(\nobj+\epsilon/2)=\wobj  + \wobj\frac{\epsilon/2}{\nobj+\epsilon/2}$; \item the particles having $\epsilon_j \ge 0$ will not be split by the S phase and will become
	$\wobj  + \wobj\frac{\nobj \epsilon_j-\epsilon/2}{\nobj+\epsilon/2}$.
	\end{itemize}
	Thus:
	\begin{equation}
		\sum_j |X_{k+1}^j - \wobj| = 
		\wobj \left[
		\sum_{\{j; \epsilon_j \le 0 \}} \frac{\epsilon/2}{\nobj+\epsilon/2}
		+
		\sum_{\{j; \epsilon_j \ge 0 \}}\frac{|\nobj \epsilon_j-\epsilon/2|}{\nobj+\epsilon/2}
		\right].
	\end{equation}
	Note that, when $\epsilon_{j_1} \ge \epsilon_{j_2} \ge 0$ replacing 
	$\epsilon_{j_2}$ by $0$ and $\epsilon_{j_1}$ 
	by $\epsilon_{j_1}+\epsilon_{j_2}$ the sum is remains constant but the error increases. Doing this to any pair of $\epsilon_j, \epsilon_{j'} \ge 0$ we obtain that the error is maximized when all of them are null except one equal to $\epsilon/2$. Thus we can write:
	\begin{equation}
		\sum_j |X_{k+1}^j - \wobj|  \le 
		\wobj	(\nobj-1) \frac{\epsilon/2}{\nobj+\epsilon/2}
		+ \wobj \frac{\nobj \epsilon/2-\epsilon/2}{\nobj+\epsilon/2} \le 
		\wobj \epsilon \frac{\nobj-1}{\nobj}.
	\end{equation}
	We conclude that, for $\epsilon < 1$: 
	$\proba [ X_{k+1} \in {\mathcal T}_{\epsilon  \frac{\nobj-1}{\nobj}} |X_k \in {\mathcal T}_{\epsilon} ] \ge 1-\epsilon$.
	Continuing this procedure, 
	$\proba \left[ \left. X_{k+\ell} \in {\mathcal T}_{\epsilon  \left( \frac{\nobj-1}{\nobj} \right)^\ell}  \right|X_k \in {\mathcal T}_{\epsilon} \right]\ge\prod_{a=1}^\ell \left[1-\epsilon\left( \frac{\nobj-1}{\nobj} \right)^a \right]$ and therefore convergence is attained with probability at least 
	$\prod_{a=1}^\infty \left[1-\epsilon\left( \frac{\nobj-1}{\nobj} \right)^a \right] \ge 1- \epsilon \nobj$ which establishes the conclusion. 
\end{proof}

$
\
$

$\star$ Then, 
in view of lemma \ref{lemma:stability},
if we prove that 
$\tau_\epsilon$ is finite almost everywhere, then the sequence $(X_l)_{l \ge 1}$  converges with probability at least $1-c_3 \epsilon$. So, 
all that remains to be proved is that 
for any $\epsilon >0$ (small enough) the stopping time $\tau_\epsilon$ is finite almost everywhere. We will fix such an $\epsilon$ from now on.
Before being able to do this, a technical result is needed.
\begin{lemma}
	If $N_0$ is at most $\nobj$ then:
	\begin{enumerate}
		\item for any $l\ge 1$: 
		\begin{equation}
			\frac{ \max_j \{ X_l^j | X_l^j \neq 0 \}   }{\min_j \{ X_l^j | X_l^j \neq 0 \} } \le 4,
		\end{equation}
		\label{item:maxmin4}
		\item for any $l\ge 0$: $N_l \le 6 \nobj$,
		\label{item:Nmax}
		\item for any $l \ge 1: \ \min \{ X_l^j | X_l^j \neq 0 \} \ge \wobj/10$.
		\label{item:wsur10}
	\end{enumerate}
	\label{lemma:finitedimension}
\end{lemma}
\begin{proof}
	\
	
	 {$\bullet$ \bf Proof of  point \ref{item:maxmin4}}: note that after a RR step the survival particles have mass exactly $\wobj$; after a S step, the resulting particles have mass in the interval $[\wobj/2, 2 \wobj]$ thus the quotient between the largest and smallest mass is at most 
	$2 \wobj / (\wobj/2)=4$. Since the renormalization step does not affect this quotient, the claim is proved for all $l\ge 1$.
	
	{$\bullet$ \bf Proof of items \ref{item:Nmax} \&  \ref{item:wsur10} }: we proceed by recurrence. 
	
	First, let us take $l=1$. Note that the RR step does not create any particles and the S step can  result in at most $2\nobj$ particles because the mass of each particle after a S step is greater than $\wobj/2$ (and total mass is conserved). Thus $N_1 \le 3 \nobj$. Moreover, the total mass ``created'' during the RR step is at most $\nobj \wobj$ thus the renormalization factor is lower bounded by $1/2$. Thus after the renormalization step each remaining particle has weight greater than $\wobj/4$.
	
	Suppose that for some $l$ both claims are true and let us prove it for $l+1$. 
	\begin{itemize}
	\item If $ \max \{ X_{l+1}^j | X_l^j \neq 0 \}  \ge \wobj$, then, using the bound 
	in item \ref{item:maxmin4} proved previously, then $X_{l+1}^k \ge \wobj/4$ for any $k \le N_{l+1}$. But since the total mass is equal to $\nobj \wobj$ it follows that there cannot be more than $4 \nobj$ particles and the claim is thus true for step $l+1$ too.
	
	\item Suppose now  $ \max \{ X_{l+1}^j | X_l^j \neq 0 \}  < \wobj$. If at the previous step all particles underwent the RR this means $N_{l+1}\le N_l \le 6\nobj$ (because RR cannot create particles) and moreover all present particles are of the same mass equal to $\nobj \wobj / N_{l+1} \ge \wobj / 6$.  
	
	\item The only alternative left is when at step $l$ there was at least a particle undergoing S, that is 
	$\max \{ X_{l}^j | X_l^j \neq 0 \}  \ge \wobj$; but, as seen above, this means that 
	$\min \{ X_{l}^j | X_l^j \neq 0 \}  \ge \wobj/4$ and therefore $N_l \le 4\nobj$ particles; the RR step cannot create particles and the S step can only result in $2\nobj$ particles (because the total mass subjected to S is at most $\nobj \wobj$ and the minimal weight after split is $\wobj/2$). Thus $N_{l+1}\le 6 \nobj$. Moreover, the total created mass by the RR step is at most $4 \wobj \nobj$ and therefore the renormalization coefficient is at least $\nobj/(\nobj+4\nobj \wobj)= 1/5$; the minimal mass before renormalization being $\wobj/2$ we conclude that the minimal mass after renormalization is $\wobj/10$ which concludes the recurrence proof.
\end{itemize}
\end{proof}

$
\
$

$\star$ Now, we can prove the last needed result:

\begin{lemma} There exists $\epsilon_0 >0$ such that 
	for any $\epsilon \in[0, \epsilon_0]$ there exists $L_\epsilon \in \N$ and $c_\epsilon < 1 $ with:
	\begin{equation}
		\forall l \ge 1: \ 
		\proba[ 
		X_{l+1} \notin \targete, ...,
		X_{l+L_\epsilon} \notin \targete | X_{l} \notin \targete ] \le c_\epsilon.
		\label{eq:tauepsilonfinite}
	\end{equation}
	Thus $\proba[\tau_\epsilon \ge k]$ is exponentially decreasing when $k\to \infty$ and in particular the
	stopping time $\tau_\epsilon$ is finite almost everywhere, i.e., 
	$\proba[\tau_\epsilon < \infty]=1$.
	\label{lemma:convergencetargete}
\end{lemma}
\begin{proof}
	We will only prove the 
	assertion \eqref{eq:tauepsilonfinite}, the rest being standard. We can assume $\epsilon <1$.

	Note that, by lemma \ref{lemma:finitedimension}, 
	with probability at least $(9/10)^{6 \nobj}$ all RR will preserve all particles, the total mass will increase and thus the renormalization factor will be smaller than unity.
	
	On the other hand, using the notation in the algorithm \ref{algo:meth_lcl}, we write any weight $X_l^j \ge \wobj$ in the form $ \wobj (I_j + R_j)$ where 
	$I_j \in \N \setminus \{0\}$ and $R_j \in [0,1[$. If $R_j/I_j \ge \epsilon/(2\nobj)$ 
	the probability of having a $I_j+1$ split is $R_j$ and thus with probability at least
	$\left(\frac{\epsilon}{2\nobj}\right)^{\nobj} \times (9/10)^{6 \nobj}$ 
	the renormalization factor is smaller than one and
	at step $l+1$ all weights larger than $\wobj$ are in the form $\wobj (1+p_j)$ with $p_j \le \frac{\epsilon}{2\nobj}$.
	
	If $X_{l+1} \in \targete$ already we are done. If not, let us denote $n^S_{l+1}$ the number 
	of   weights larger than $\wobj$, denoted from now one $\wobj (1+p_j)$, $j=1,..., n^S_{l+1}$. Denote also $n^R_{l+1}$ the number 
	of   weights strictly smaller than $\wobj$. By the conservation of mass, the sum of these latter weights is $$\left[(\nobj-n^S_{l+1}) -\sum_j p_j \right]\wobj \ge 
	\left[(\nobj-n^S_{l+1}) -\epsilon/2 \right]\wobj > 
	(\nobj-n^S_{l+1}-1)\wobj.$$
	Since all weights in this set are inferior to $\wobj$ we deduce that $n^R_{l+1} > \nobj-n^S_{l+1}-1$.
	
	 If $n^R_{l+1} = \nobj-n^S_{l+1}$ then $X_{l+1} \in \targete$ (situation already discussed).
	  
	Thus we can assume $n^R_{l+1} \ge \nobj-n^S_{l+1}+1$. 
	
	\noindent Consider now the $(\nobj-n^S_{l+1}+1)$-th weight counting down, in decreasing magnitude, from the largest one, denote it $y$. Since the sum of the $\nobj-n^S_{l+1}+1$ largest weights is smaller than $\left[(\nobj-n^S_{l+1})\right]\wobj$ it follows that 
	$y \le (\nobj-n^S_{l+1}) / (\nobj-n^S_{l+1}+1) \le \nobj / (\nobj+1)$. That means that with probability larger than 
	$ \left[ \frac{1}{\nobj+1}\right]^{6\nobj} \times 
	\left[ \frac{1}{10}\right]^{6\nobj}$ at the step $l+2$ we can keep, after RR, exactly 
	the $\nobj-n^S_{l+1}$ largest weights and eliminate the others. Doing this the total mass will increase, with at most $\epsilon/2$, therefore the renormalization coefficient will be smaller than one and larger than 
	$\frac{\nobj}{\nobj + \epsilon/2}$. Note that all particles $\wobj(1+p_i)$ are kept unchanged by the S step with probability $1-\epsilon/2$. Therefore, with the strictly positive probability 
	$(1-\epsilon/2) \left[ \frac{1}{\nobj+1}\right]^{6\nobj} \times 
	\left[ \frac{1}{10}\right]^{6\nobj}$
	(depending only on $\epsilon$ and $\nobj$) we will have exactly $\nobj$ particles, and those over $\wobj$ have a distance to $\wobj$ at most $\epsilon/(2\nobj)$, which means $X_{l+2} \in \targete$. The lemma is proved with $L_\epsilon=2$.
\end{proof}

$
\
$

$\star$ We can now finish the proof of the proposition \ref{prop:cv_conservative}: we showed (lemma \ref{lemma:returntonobj}) that, after a finite time $\tau$,  $X_l$ will have at most $\nobj$ non-null particles and, by lemma \ref{lemma:finitedimension} is will have at most $\nmax= 6 \nobj$ non-null particles forever from that point on. From lemma \ref{lemma:convergencetargete} we conclude that after a finite time $\tau_\epsilon$ the state $X_{\tau_\epsilon}$ will be in $\targete$ and from lemma \ref{lemma:stability} it will converge with probability $1-c_3 \epsilon $. Since this is true for any $\epsilon$ (small enough) then we have convergence everywhere, hence the conclusion.

}

\newpage

\bibliographystyle{plain}



\end{document}